\documentclass[10pt,reqno]{amsart}
\usepackage{amsmath,amsopn,amssymb,amsthm,multicol}
\usepackage{hyperref}
\usepackage{color}
\usepackage{graphicx}

\DeclareMathOperator{\Ad}{Ad}
\DeclareMathOperator{\ad}{ad}

\DeclareMathOperator{\Mark}{Mrk}
\DeclareMathOperator{\Span}{span}

\newcommand{\fr}{\mathfrak}
\newcommand{\al}{\alpha}
\newcommand{\be}{\beta}

\newcommand{\bb}{\mathbb}
\DeclareMathOperator{\rnk}{rk}

\DeclareMathOperator{\SU}{SU}
\DeclareMathOperator{\U}{U}
\DeclareMathOperator{\G}{G}
\DeclareMathOperator{\F}{F}
\DeclareMathOperator{\E}{E}
\newcommand{\thickline}{\noalign{\hrule height 1pt}}

 \newtheorem{lemma} {Lemma} [section]
\newtheorem{theorem}[lemma]{Theorem} 
\newtheorem{remark}[lemma] {Remark} 
\newtheorem{prop} [lemma]{Proposition}  
\newtheorem{definition}[lemma] {Definition}

\begin{document}

\title{Equigeodesics on generalized flag manifolds with $\G_2$-type $\fr{t}$-roots} 
\author{Marina Statha}
\address{University of Patras, Department of Mathematics, GR-26500 Rion, Greece}
\email{statha@master.math.upatras.gr} 
\medskip

\begin{abstract}
We study homogeneous curves in generalized flag manifolds $G/K$ with $\G_2$-type $\fr{t}$-roots, which are geodesics with respect to each $G$-invariant metric on $G/K$.  These curves are called equigeodesics.  The tangent space of such flag manifolds splits into six isotropy summands, which are in  one-to-one correspondence with $\fr{t}$-roots.  Also, these spaces are a generalization of the exceptional full flag manifold $\G_2/T$.  We give a characterization for structural equigeodesics for flag manifolds with $\G_2$-type $\fr{t}$-roots, and we give for each such flag manifold, a list of subspaces in which the vectors are structural equigeodesic vectors.

\medskip
\noindent 2010 {\it Mathematics Subject Classification.} Primary 53C25; Secondary 53C30, 13P10, 65H10, 68W30.

\medskip
\noindent {\it Keywords}:  Homogeneous space, flag manifold, equigeodesic vector, $\fr{t}$-root.
\end{abstract}

\maketitle
 

\section{Introduction}
\markboth{Marina Statha}{Equigeodesics on generalized flag manifolds with $\G_2$-type $\fr{t}$-roots}

Let $(G/K, g)$ be a Riemannian homogeneous space.  A geodesic $\gamma(t)$ through the origin $o=eK$ is called a {\it homogeneous geodesic} if it is an orbit of a one-parameter subgroup of $G$, i.e.,
$
\gamma(t)=(\exp tX)\cdot o,
$
where $X$ is a non-zero vector in the Lie algebra $\fr{g}$ of $G$.  If all geodesics on $G/K$ are homogeneous geodesics the homogeneous space is called a {\it g.o. manifold} (from ``geodesic ordit").  The terminology was introduced by O. Kowalski and L. Vanhecke in \cite{KoVa}, who initiated a systematic study of such spaces.  Examples of such spaces are the symmetric spaces, the weakly symmetric spaces and the naturally reductive spaces.  

In \cite{CGN} the authors studied homogeneous curves on generalized flag manifolds that are geodesics with respect to each invariant metric on the flag manifold.  These curves are called {\it equigeodesics}.  Since the infinitesimal generator of the one parameter subgroup is an element of the Lie algebra of $G$, it is natural to characterize the equigeodesics in terms of their infinitesimal generator.  This allows us to use a Lie theoretical approach for the study of homogeneous geodesics on flag manifolds.  The infinitesimal generator of an equigeodesics is called {\it equigeodesic vector}.  An algebraic characterization of equigeodesic vectors on generalized flag manifolds is given in \cite{CGN}.  

Recall that a {\it generalized flag manifold} is a homogeneous space $G/K$ where $G$ is a compact, semisimple Lie group and $K$ is the centralizer of a torus in $G$.  Actually a vector is equigeodesic if and only if it is a solution of an algebraic system of equations whose variables are the components of the vector.  However, there exist some subspaces of the tangent space $\fr{m}\cong T_{o}(G/K)$ of the flag manifold $G/K$, all of whose elements are equigeodesic vectors.  The existence of such subspaces depends on the geometric structure of the $G/K$.  These equigeodesic vectors are called {\it structural equigeodesic}.  The authors in \cite{GrNe} have provided a version of the previously formula for equigeodesic vectors on generalized flag manifolds with two isotropy summands.  Later in \cite{WaZh} the authors gave a general formula for finding equigeodesic vectors on generalized flag manifolds with second Betti number equal to one (that is flag manifolds which are determined by painting one black node their Dynkin diagram).  More precisely, they found families of subspaces in which all vectors are structural equigeodesic vectors, on generalized flag manifolds associated to exceptional Lie groups $\F_4, \E_6$ and $\E_7$ with three isotropy summands, that is $\F_4/(\U(2)\times\SU(3))$, $\E_6/(\U(2)\times\SU(3)\times\SU(3))$ and $\E_7/(\U(3)\times\SU(5))$. 

In the present article we study equigeodesics on generalized flag manifolds with $\G_2$-type $\fr{t}$-roots.  In particular, for such type of flag manifolds we describe the families of subspaces in which all elements are structural equigeodesic vectors.  We know from \cite{ArChSa2} that generalized flag manifolds $G/K$ with $\G_2$-type $\fr{t}$-roots have six isotropy summands and correspond to painted Dynkin diagrams with two black nodes with Dynkin marks $2$ and $3$.  In particular, these are the generalized flag manifolds $\F_4(\alpha_3, \al_4), \E_6(\al_3, \al_6)$, $\E_7(\al_5, \al_6)$ and $\E_8(\al_1, \al_2)$.  Here with $G(\al_{i_{0}}, \al_{j_{0}})$ we denote the flag manifold $M=G/K$ where we have painted two black nodes on the Dynkin diagram of $G$.

In Theorems \ref{BasicTheorem} and \ref{BasicTheorem1} of the present paper we provide a method to obtain structural equigeodesic vectors (cf. Propositions \ref{prop4.1} and \ref{prop4.2}).  For $\F_4(\alpha_3, \al_4)$ and $\E_6(\al_3, \al_6)$ we find all subspaces in which the vectors are structural equigeodesics.  These are described in Tables 3,4 and 5.  For the flag manifold $\E_7(\al_5, \al_6)$ we describe all the roots that satisfy Theorem \ref{BasicTheorem} and therefore we can describe by simple calculation all root spaces whose vectors are structural equigeodesics (we give some of them).  Finally, for the flag manifold $\E_8(\al_1, \al_2)$ we give some of the roots that satisfy Theorem \ref{BasicTheorem1}. 
In conclusion we have the following:
\begin{theorem}
The generalized flag manifolds $\F_4/(\U(3)\times\U(1))$, $\E_6/(\U(3)\times\U(3))$, $\E_7/(\U(6)\times\U(1))$ and $\E_8/(\E_6\times\U(1)\times\U(1))$ admit non trivial structural equigeodesic vectors.
\end{theorem}

\section{Generalized Flag Manifolds}

\subsection{Description of flag manifolds in terms of painted Dynkin diagrams} 
Let $\fr{g}$ and $\fr{k}$ be the Lie algebras of $G$ and $K$ respectively and $\fr{g}^{\bb{C}}$, $\fr{k}^{\bb{C}}$ be their complexifications.  We choose a maximal torus  $T$ in $G$ and let $\fr{h}$ be the Lie algebra of $T$.  Then the complexification $\fr{h}^{\bb{C}}$ is a Cartan subalgebra of $\fr{g}^{\bb{C}}$.  Let $R\subset(\fr{h}^{\bb{C}})^{*}$   be the root system of $\fr{g}^{\mathbb{C}}$ relative to the Cartan subalgebra $\fr{h}^{\mathbb{C}}$  and consider the  root space decomposition $\fr{g}^{\mathbb{C}}=\fr{h}^{\mathbb{C}}\oplus\sum_{\al\in R}\fr{g}_{\al}^{\mathbb{C}}$, where $\fr{g}_{\al}^{\bb{C}}=\{X\in\fr{g}^{\bb{C}} : \ad(H)X=\al(H)X,   \ \mbox{for all} \   H\in\fr{h}^{\bb{C}}\}$  denotes the root space associated to a root $\al$.  Assume that $\fr{g}^{\bb{C}}$ is semisimple, so the Killing form $B$ of $\fr{g}^{\bb{C}}$ is non degenerate, and we establish a natural isomorphism between $\fr{h}^{\bb{C}}$ and the dual space $(\fr{h}^{\bb{C}})^{*}$ as follows: for every $\al\in(\fr{h}^{\bb{C}})^{*}$ we define $H_{\al}\in\fr{h}^{\bb{C}}$ by the equation $B(H, H_{\al}) = \al(H)$, for all $H\in\fr{h}^{\bb{C}}$.  We take a Weyl basis  $E_{\al}\in\fr{g}_{\al}^{\mathbb{C}}$ $(\al\in R)$  with $B(E_{\al}, E_{-\al})=-1$  and  $[E_{\al}, E_{-\al}]=-H_{\al}$. Then   $\fr{g}_{\al}^{\mathbb{C}}=\mathbb{C}E_{\al}$  and 
\begin{equation}\label{LieBracket}
 [E_{\al}, E_{\be}]=
 \left\{
\begin{array}{ll}
  N_{\al, \be}E_{\al+\be}  & \mbox{if} \ \ \al, \be,\al+\be\in R \\
  0 & \mbox{if} \ \ \al, \be\in R,  \al+\be\notin R,
\end{array} \right.
\end{equation}
where the structure constants $N_{\al, \be}\in\bb{R}$ are such that $N_{\al, \be}=0$  if $\al, \be \in R$, $\al+\be\notin R$, and $N_{\al, \be}=-N_{\be, \al}$,  $ N_{\al, \be}=N_{-\al, -\be}\in \mathbb{R}$ if 
$\al, \be, \al+\be\in R$.  It is clear  that $N_{\al, \be}\neq 0$, if $\al, \be, \al+\be\in R$  and so relation (\ref{LieBracket}) implies  that   $[\fr{g}_{\al}^{\bb{C}}, \fr{g}_{\be}^{\bb{C}}]=\fr{g}_{\al+\be}^{\bb{C}}$.
Choose a  basis $\Pi=\{\al_{1}, \ldots, \al_{\ell}\}$ $(\dim\fr{h}^{\mathbb{C}}=\ell)$ of simple roots for $R$, and let $R^{+}$ be a choise of positive roots.  Set $A_{\al}=E_{\al}+E_{-\al}$ and  $B_{\al}=\sqrt{-1}(E_{\al}-E_{-\al})$, where $\al\in R^{+}$. Then the real subalgebra  $\fr{g}$ is  given by
\begin{equation}\label{realg}
\fr{g}=\fr{h}\oplus \sum_{\al\in R^{+}}(\mathbb{R}A_{\al}+\mathbb{R}B_{\al}) = \fr{h}\oplus\sum_{\al\in R^{+}}\fr{U}_{\al}.
\end{equation}
 Note that $\fr{g}$, as a real form of $\fr{g}^{\bb{C}}$ is the fixed point set of the conjugation   $\tau : \fr{g}^{\bb{C}}\to\fr{g}^{\bb{C}}$, which without loss of generality   can be assumed to be such  that $\tau(E_{\al})=E_{-\al}$.

Since $\fr{h}^{\bb{C}}\subset\fr{k}^{\bb{C}}\subset\fr{g}^{\bb{C}}$,  there is a closed subsystem $R_{K}$ of $R$ such that 
$\fr{k}^{\bb{C}}=\fr{h}^{\bb{C}}\oplus\sum_{\al\in R_{K}}\fr{g}_{\al}^{\mathbb{C}}$.  In particular, we can always  find  a subset $\Pi_{K}\subset\Pi$   such that  
$
R_{K}=R\cap\left\langle\Pi_{K} \right\rangle=\{\be\in R : \be=\sum_{\al_{i}\in\Pi_{K}}k_{i}\al_{i}, \ k_{i}\in\bb{Z}\},
$
where    $\left\langle\Pi_{K} \right\rangle$  is the space of roots generated by $\Pi_{K}$  with integer coefficients. The complex Lie algebra $\fr{k}^{\mathbb{C}}$  is a maximal reductive subalgebra of $\fr{g}^{\mathbb{C}}$  
and thus it  admits the decomposition $\fr{k}^{\mathbb{C}}=\fr{z}(\fr{k}^{\mathbb{C}})\oplus\fr{k}^{\mathbb{C}}_{ss}$,   
where $\fr{z}(\fr{k}^{\mathbb{C}})$ is  the  center of $\fr{k}^{\mathbb{C}}$ and  $\fr{k}^{\mathbb{C}}_{ss}=[\fr{k}^{\mathbb{C}}, \fr{k}^{\mathbb{C}}]$ 
is its  semisimple part.  Note that $\fr{k}^{\bb{C}}_{ss}$ is given by
$\fr{k}^{\bb{C}}_{ss}=\fr{h}'\oplus\sum_{\al\in R_{K}}\fr{g}_{\al}^{\bb{C}},$
where $\fr{h}'=\sum_{\al\in\Pi_{K}}\bb{C}H_{\al}\subset\fr{h}^{\bb{C}}$ is a Cartan subalgebra of $\fr{k}^{\mathbb{C}}_{ss}$.  In fact, $R_{K}$ is  the root system of the semisimple part  $\fr{k}^{\mathbb{C}}_{ss}$ and  $\Pi_{K}$ is a corresponding basis.  Thus we easily conclude that $\dim_{\bb{C}}\fr{h}'={\rm card}\,{\Pi_{K}}$, where  ${\rm card}\,{\Pi_{K}}$ denotes  the cardinality of the set $\Pi_{K}$.
Let $K$ be the connected Lie subgroup of $G$ generated by $\fr{k}=\fr{k}^{\bb{C}}\cap \fr{g}$.  Then the homogeneous manifold $M=G/K$ is a flag manifold, and any flag manifold is defined in this way, i.e. by the choise of a triple $(\fr{g}^{\bb{c}}, \Pi, \Pi_{K})$.

Set  $\Pi_{M}=\Pi\backslash \Pi_{K}$, and  $R_{M}=R\backslash R_{K}$, such that $\Pi=\Pi_{K}\cup \Pi_{M}$, and $R=R_{K}\cup R_{M}$, respectively.  Roots in $R_{M}$ are called  {\it complementary roots}, and they play an important role in the geometry of $M=G/K$.  
For example, let $\fr{m}$ the orthogonal complement of $\fr{k}$ in $\fr{g}$ with respect to $B$.  Then we have $[\fr{k}, \fr{m}]\subset\fr{m}$ where $\fr{m}\cong T_{o}(G/K)$.  We set $R_{M}^{+} = R^{+}\backslash R_{K}^{+}$ where $R_{K}^{+}$ is the system of positive roots of $\fr{k}^{\bb{C}}$ $(R_{K}^{+}\subset R^{+})$.  Then 
\begin{equation}\label{tangent}
\fr{m} = \sum_{\al\in R_{M}^{+}}(\bb{R}A_{\al} + \bb{R}B_{\al}).
\end{equation} 
The complexification is given as $\fr{m}^{\mathbb{C}}=\sum_{\al\in R_{M}}\mathbb{C}E_{\al}$, and the set $\{E_{\al} : \al\in R_{M}\}$ is a basis of $\fr{m}^{\bb{C}}$.

Now if we assume that $\Pi_{M}=\Pi\backslash\Pi_{K} = \{\al_{i_1},\ldots,\al_{i_{r}}\}$, where $1\leq i_{1}<\cdots<i_{r}\leq\ell$ we set, for some integers $j_1, \dots, j_r$ with $(j_1, \dots, j_r) \neq (0, \dots, 0)$ 
\begin{equation}\label{summands}
R^{\fr{m}}( j_1, \dots, j_r ) = \left\{\, \sum_{j=1}^{l} m_j \alpha_j \in R^{+} \  :  \  m_{i_1} = j_1, \dots, m_{i_r} = j_r \, \right\}\subset R^{+}. 
\end{equation}
Note that $R_{M}^{+} =  R^+ \backslash R_{K}^{+}= \bigcup_{j_1, \dots, \, j_r} R^{\fr{m}}( j_1, \dots, j_r )$.  For any $R^{\fr{m}}( j_1, \dots, j_r ) \neq \emptyset$, we define an $\Ad(K)$-invariant subspace ${\frak m}( j_1, \dots, j_r )$ of $\frak g$ by 
\begin{equation}\label{summands2}
{\frak m}( j_1, \dots,  j_r ) = \sum_{\alpha \in  R^{\fr{m}}( j_1, \dots, \, j_r )}
 \left\{ {\mathbb R}A_{\al} + {\mathbb R}B_{\al} \right\}. 
\end{equation}
Then we have a decomposition of $\frak m$ into mutually non equivalent irreducible $\Ad(K)$-modules ${\frak m}( j_1, \dots,  j_r )$:
$$
\frak m = \sum_{j_1, \dots, \, j_r} {\frak m}( j_1, \dots,  j_r ). 
$$

We conclude that all information contained in  $\Pi=\Pi_{K}\cup\Pi_{M}$ can be presented graphically by the  painted Dynkin diagram of $M=G/K$.
\begin{definition}\label{pdd}
Let $\Gamma=\Gamma(\Pi)$ be the Dynkin diagram of the fundamental system $\Pi$.  By painting in  black the nodes 
of $\Gamma$  corresponding  to  $\Pi_{M}$, we obtain the painted Dynkin diagram of the flag manifold $G/K$. In this diagram
the subsystem $\Pi_{K}$ is determined as the subdiagram of white roots. 
\end{definition}
  
Conversely, given a painted Dynkin diagram, in order to obtain the corresponding flag manifold $M=G/K$ we are working as follows: We define $G$ as the unique simply connected Lie group  corresponding to the  underlying Dynkin diagram $\Gamma=\Gamma(\Pi)$. The connected Lie subgroup $K\subset G$ is defined by using  the additional information $\Pi=\Pi_{K}\cup\Pi_{M}$ encoded into the painted Dynkin diagram.  The semisimple part of $K$ is obtained from  the (not  necessarily connected) subdiagram of white roots, and   each black root, i.e.   each  root in $\Pi_{M}$,  gives rise to one $U(1)$-summand.  Thus the painted Dynkin diagram determines  the  isotropy subgroup $K$ and the space $M=G/K$ completely.   By using certain rules to determine whether  different painted Dynkin diagrams define isomorphic flag manifolds (see \cite{AlAr}),  one  can obtain all flag manifolds $G/K$ of  a compact  simple Lie group $G$.  

From now on we denote the flag manifold $M=G/K$ with $G\in\{B_{\ell}, C_{\ell}, D_{\ell}, \F_4, \E_6, \E_7, \E_8\}$, by $G(\al_{i_{0}})$ if we have painted one node of $\Gamma(\Pi)$, that is $\Pi_{M} = \Pi\backslash \Pi_{K}=\{\al_{i_{0}}\}$ and by $G(\al_{i_{0}}, \al_{j_{0}})$ if we have painted two nodes of $\Gamma(\Pi)$, that is $\Pi_{M}=\Pi\backslash\Pi_{K}=\{\al_{i_{0}}, \al_{j_{0}}\}$.

We close this subsection with the next lemma which gives us some information about the Lie algebra structure of $\fr{g}$.

\begin{lemma}\label{basiclemma}
The Lie brackets among the elements of the basis $\{A_{\al}, B_{\al}, \sqrt{-1}H_{\beta}: \al\in R^{+}\ \mbox{and}\ \beta\in \Pi\}$ of $\fr{g}$ are given as follows:
\begin{itemize}
\item[] $[A_{\al}, A_{\beta}] = N_{\al,\beta}\, A_{\al+\beta} + N_{-\al,\beta}\, A_{\al-\beta}$, \quad \quad $[\sqrt{-1}H_{\al}, A_{\beta}] = \beta(H_{\al})B_{\beta}$

\item[] $[B_{\al}, B_{\beta}] = -N_{\al, \beta}\, A_{\al+\beta} - N_{\al, -\beta}\, A_{\al-\beta}$, \quad $[\sqrt{-1}H_{\al}, B_{\beta}]=-\beta(H_{\al})A_{\beta}$

\item[] $[A_{\al}, B_{\beta}] = N_{\al,\beta}\, B_{\al+\beta} + N_{\al,-\beta}\, B_{\al-\beta}$, \quad \quad $[A_{\al}, B_{\al}]=2\sqrt{-1}H_{\al}$,
\end{itemize}
where $\al+\beta$, $\al-\beta$ are roots.
\end{lemma} 
\begin{proof}
We will prove three of the above relations and the others can be obtained by a similar method.  For the first we have:
\begin{eqnarray*}
[A_{\al}, A_{\beta}] &=& [E_{\al}+E_{-\al}, E_{\beta}+E_{-\beta}] = N_{\al,\beta}E_{\al+\beta} + N_{\al,-\beta}E_{\al-\beta} + N_{-\al,\beta}E_{-\al+\beta} + N_{-\al,-\beta}E_{-\al-\beta} \\ 
&=& N_{\al,\beta}E_{\al+\beta} + N_{-\al, \beta}E_{\al-\beta} + N_{-\al,\beta}E_{-\al+\beta} + N_{\al, \beta}E_{-(\al+\beta)} \\ 
&=& N_{\al,\beta} A_{\al+\beta} + N_{-\al, \beta}A_{\al-\beta}.
\end{eqnarray*}
 
For the second we have:
\begin{eqnarray*}
[\sqrt{-1}H_{\al}, A_{\beta}] &=& [\sqrt{-1}H_{\al}, E_{\beta}+E_{-\beta}] = \sqrt{-1}[H_{\al}, E_{\beta}] + \sqrt{-1}[H_{\al}, E_{-\beta}] \\
&=& \sqrt{-1} \beta(H_{\al})E_{\beta} - \sqrt{-1}\beta(H_{\al})E_{-\beta}\\
&=& \beta(H_{\al})B_{\beta}.
\end{eqnarray*}
Finally, we prove the last relation:
\begin{eqnarray*}
[A_{\al}, B_{\al}] &=& [E_{\al}+E_{-\al}, \sqrt{-1}(E_{\al}-E_{-\al})] = -\sqrt{-1}[E_{\al}, E_{-\al}] + \sqrt{-1}[E_{-\al}, E_{\al}] \\
&=& \sqrt{-1}H_{\al} + \sqrt{-1}H_{\al} = 2\sqrt{-1}H_{\al}.
\end{eqnarray*}
\end{proof}

\subsection{$\fr{t}$-roots and isotropy summands}
We study the isotropy representation of a generalized flag manifold $M=G/K$ of a compact simple Lie group $G$  in terms of $\fr{t}$-roots.  In order to realise the decomposition of $\fr{m}$ into irreducible $\Ad(K)$-modules we use the center $\fr{t}$ of the real Lie algebra $\fr{k}$.  For simplicity, we fix  a system of simple roots $\Pi=\{\al_1,\ldots, \al_r, \phi_1, \ldots, \phi_k\}$ of $R$,  such that $r+k=\ell=\rnk\fr{g}^{\bb{C}}$  and we assume that $\Pi_{K}=\{\phi_1, \ldots, \phi_k\}$ is a basis of the root system $R_{K}$ of $K$ so $\Pi_{M}=\Pi\backslash \Pi_{K}=\{\al_{1}, \ldots, \al_{r}\}$.  Let $\Lambda_{1}, \ldots, \Lambda_{r}$ be  the fundamental weights  corresponding to the simple roots of $\Pi_{M}$, i.e. the linear forms defined by 
$\frac{2(\Lambda_{i}, \al_{j})}{(\al_{j}, \al_{j})}=\delta_{ij}, (\Lambda_{j}, \phi_{i})=0,$
where $(\al, \be)$ denotes the inner product on $(\fr{h}^{\mathbb{C}})^{*}$ given by $(\al, \be)=(H_{\al}, H_{\be})$, for all $\al, \be\in(\fr{h}^{\mathbb{C}})^{*}$.  Then the $\{\Lambda_{i} : 1\leq i\leq r\}$ is a basis of the dual space $\fr{t}^{*}$ of $\fr{t}$,  $\fr{t}^{*}=\sum_{i=1}^{r}\mathbb{R}\Lambda_{i}$ and $\dim\fr{t}^*=\dim \fr{t}=r$.

Consider now the linear restriction map  $\kappa : \fr{h}^{*}\to \fr{t}^{*}$ defined by $\kappa(\al)=\al|_{\fr{t}}$,
and set  $R_{\fr{t}} = \kappa(R)=\kappa(R_{M})$.  
\begin{definition}
The elements of $R_{\fr{t}}$ are called {\it $\fr{t}$-roots}.  
\end{definition}
The set $R_{\fr{t}}$ is not in general a root system.  An element $Y\in\fr{t}$  is called {\it regular}  if any $\fr{t}$-root $\kappa(\al) =\xi$ $(\al\in R_{M})$ takes non zero value at $Y$, i.e. $\xi(Y)\neq 0$.   A regular element defines an ordering
in $\fr{t}^*$  and thus we obtain the splitting  $R_{\fr{t}}=R_{\fr{t}}^{+}\cup  R_{\fr{t}}^{-}$,  where $R_{\fr{t}}^{+}=\{\xi\in R_{\fr{t}} : \xi(Y)>0\}$ and $R_{\fr{t}}^{-}=\{\xi\in R_{\fr{t}} : \xi(Y)<0\}$.  The $\fr{t}$-toots $\xi\in R_{\fr{t}}^{+}$ (resp. $\xi\in R_{\fr{t}}^{-}$) will be called {\it positive} (resp. {\it negative}). Since $R_{\fr{t}}=\kappa(R_{M})$ it follows that $R_{\fr{t}}^{+}=\kappa(R_{M}^{+})$, and since $R_{M}^{-}=-R_{M}^{+}=\{-\al : \al\in R_{M}^{+}\}$,  it is  $R_{\fr{t}}^{-}=\kappa(R_{M}^{-})$. 

\begin{definition}
A $\fr{t}$-root is called simple if it is not a sum of two positive $\fr{t}$-roots.  
\end{definition}
The set of all simple $\fr{t}$-roots is denoted as $\Pi_{\fr{t}}$ and is a basis of $\fr{t}^*$,  in the sense that any $\fr{t}$-root can be written  as a linear combination of its elements with integer coefficients of the same sign. We will call the set $\Pi_{\fr{t}}$ as a {\it $\fr{t}$-basis}.
\begin{prop}\label{tbasis} \textnormal{(\cite{ArCh2})}
A $\fr{t}$-basis $\Pi_{\fr{t}}$  is obtained by restricting  the roots of $\Pi_{M}=\Pi\backslash 
\Pi_{K}$ to $\fr{t}$, that is
$
\Pi_{\fr{t}}=\{\kappa(\al_i)=\overline{\al}_{i}=\al_{i}|_{\fr{t}} : \al_{i}\in \Pi_{M}\}.
$
\end{prop} 

As we saw the flag manifolds $G/K$ are determined by pairs $(\fr{g}, \Pi, \Pi_K)$.  The number of $\ad(\fr{k})$-submodules of $\fr{m}\cong T_{o}(G/K)$ correspond to the {\it Dynkin mark} of the simple root we paint black on the Dynkin diagram.  We recall the following definition

\begin{definition}\label{mark}
The Dynkin mark of a simple root $\al_{i}\in\Pi$  $(i=1, \ldots, \ell)$, is  the  positive integer $m_{i}$ in 
the expression of the highest root $\widetilde{\al}=\sum_{k=1}^{\ell}m_{k}\al_{k}$ in terms of simple roots.  We will denote by $\Mark$ the function $\Mark : \Pi\to\bb{Z}^{+}$  with $\Mark(\al_{i})=m_{i}$.
\end{definition}

By using the Proposition \ref{tbasis} we give a useful method to find the positive $\fr{t}$-root $R_{\fr{t}}^+$ of $\Pi_M = \{\al_i\in \Pi : \Mark(\al_i) = m_i\}$.  The $\fr{t}$-basis is $\Pi_{\fr{t}} = \{\overline{\al}_i\}$, where $\kappa(\al_i) = \overline{\al}_i=\al_i|_{\fr{t}}$ and $\fr{t}^{*}=\bb{R}\overline{\al}_i$.  We fix a positive root $\al=\sum_{j=1}^{\ell}k_j\al_j\in R^{+}$, with $k_j\leq m_j$  for all $j=1,\ldots, \ell$.  Then by the fact that $\kappa(R_K) = 0$ we have that for all $\al\in R_{M}$, $\kappa(\al) = \al|_{\fr{t}} = k_i\overline{\al}_i$ with $1\leq k_i\leq m_i$.  Hence $R_{\fr{t}}^{+} = \{k_i\overline{\al}_i : 1\leq k_i\leq m_i\}$ $=\{\overline{\al}_i, 2\overline{\al}_i,\ldots,m_i\overline{\al}_i\}$ and ${\rm card} R_{\fr{t}}^+ = m_i$.  In case where $\Pi_M=\{\al_i, \al_j : \Mark(\al_i) = m_i, \Mark(\al_j)=m_j\}$ then $\Pi_{\fr{t}} = \{\overline{\al}_i, \overline{\al}_j : i<j\}$, where $\kappa(\al_i) = \overline{\al}_i=\al_i|_{\fr{t}}$, $\kappa(\al_j)=\overline{\al}_j=\al_j|_{\fr{t}}$ and $\fr{t}^{*} =\Span\{\overline{\al}_i, \overline{\al}_j\}$.  Then for $\al\in R^{+}_{M}$ we have $\kappa(\al)=\al\big|_{\fr{t}}=k_{i}\overline{\al}_{i}+k_{j}\overline{\al}_{j}$ where $0\leq k_{i}\leq m_{i},\ 0\leq k_{j}\leq m_{j}$
the coefficients $k_i, k_j$ can not be simultaneously zero, so it is obvious that ${\rm card} R_{\fr{t}}^{+}\geq 3$.  For generalized flag manifold $G/K$ with $b_2(G/K)\geq 3$ there are more than five $\fr{t}$-roots (\cite{ArChSa1}). 

A fundamental result about $\fr{t}$-root is the following:

\begin{prop}\textnormal{(\cite{AlPe})} \label{isotropy} 
There exists a one-to-one correspondence between $\fr{t}$-roots $\xi$ and irreducible $\ad(\fr{k}^{\mathbb{C}})$-submodules $\fr{m}_{\xi}$\footnote{We mean that $[\fr{k}^{\bb{C}}, \fr{m}_{\xi}]\subset\fr{m}_{\xi}$ for all $\xi\in R_{\fr{t}}$.} of the isotropy representation of $\fr{m}^{\mathbb{C}}$, which is given by   
$$
R_{\fr{t}}\ni\xi \ \leftrightarrow \ \fr{m}_{\xi} =\sum_{\al\in R_{M}: \kappa(\al)=\xi}\mathbb{C}E_{\al}.
$$
Thus $\fr{m}^{\mathbb{C}} = \bigoplus_{\xi\in R_{\fr{t}}} \fr{m}_{\xi}.$  Moreover, these submodules are non equivalent as $\ad(\fr{k}^{\mathbb{C}})$-modules.		 		   
\end{prop}

In order to obtain a decomposition of the real $\Ad(K)$-module $\fr{m}$ in terms of $\fr{t}$-roots, we use the complex conjugation $\tau$ of $\fr{g}^{\bb{C}}$ with respect to $\fr{g}$ (note that $\tau$ interchanges $\fr{g}_{\al}^{\bb{C}}$ and $\fr{g}_{-\al}^{\bb{C}}$).  Moreover, for a complex subspace $V$ of $\fr{g}^{\bb{C}}$ we denote by $V^{\tau}$ the set of all fixed points of $\tau$.  Then, we can write
\begin{equation}\label{complextan}
\fr{m} = \bigoplus_{\xi\in R_{\fr{t}}^{+}}(\fr{m}_{\xi}\oplus\fr{m}_{-\xi})^{\tau}.
\end{equation}   

Let us assume for simplicity that $R^{+}_{\fr{t}} = \{\xi_1,\dots, \xi_{s}\}$.  In this case Proposition \ref{isotropy} and relations (\ref{tangent}), (\ref{complextan}) imply that the real irreducible $\ad(\fr{k})$-submodule $\fr{m}_i=(\fr{m}_{\xi_{i}}\oplus\fr{m}_{-\xi_{i}})^{\tau}$ $(1\leq i\leq s)$ which corresponds to a positive $\fr{t}$-root $\xi_{i}$, is necessarily of the form
\begin{equation}\label{sxesi}
\fr{m}_i = \sum_{\{\al\in R_{M}^{+} : \kappa(\al) = \xi_i\}} \{\bb{R}A_{\al} + \bb{R}B_{\al}\}.
\end{equation}

By summarizing, we have the following proposition

\begin{prop}\textnormal{(\cite{ArChSa1})}\label{onetone}
Let $M=G/K$ be a generalized flag manifold defined by a subset $\Pi_{K}\subset\Pi$ such that $\Pi_{M}=\Pi\backslash\Pi_{K}=\{\al_{i_{1}},\dots,\al_{i_{r}}\}$ with $1\leq i_1<\cdots <i_{r}\leq\ell$.  Assume that $\fr{g}=\fr{k}\oplus\fr{m}$ is a $B$-orthogonal reductive decomposition.  Then
\begin{enumerate}
\item There exists a natural one-to-one correspondence between elements of the set $R^{\fr{m}}(j_1,\dots,j_r)$ and the set of positive $\fr{t}$-roots $R_{\fr{t}}^{+} = \{\xi_1,\dots,\xi_s\}$.  Therefore, there is a decomposition of $\fr{m}$ into $s$ mutually non-equivalent irreducible $\Ad(K)$-modules
$
\fr{m} = \sum_{\xi\in R_{\fr{t}}^{+}}(\fr{m}_{\xi}\oplus\fr{m}_{-\xi})^{\tau} = \sum_{i=1}^{s}(\fr{m}_{\xi_{i}}\oplus\fr{m}_{-\xi_{i}})^{\tau} = \sum_{j_1,\dots, j_r}\fr{m}(j_1,\dots,j_r)
$
\item The dimensions of the real $\Ad(K)$-modules $\fr{m}_i$ $(i=1,\ldots, s)$ corresponding to the $\fr{t}$-root $\xi_{i}\in R_{\fr{t}}^{+}$ are given by $\dim_{\bb{R}}\fr{m}_i = 2\, {\rm card}\,\{\al\in R_{M}^{+} : \kappa(\al) = \xi_{i}\}$ $=2\, {\rm card}R^{\fr{m}}\,(j_1,\dots j_r)$, for appropriate positive integers $j_1,\dots j_r$.
\item Any $G$-invariant Riemannian metric $g$ on $G/K$ is given by
$$
g = \sum_{\xi\in R^{+}_{\fr{t}}}x_{\xi}B|_{(\fr{m}_{\xi}\oplus\fr{m}_{-\xi})^{\tau}} = \sum_{i=1}^{s}x_{\xi_{i}}B|_{(\fr{m}_{\xi_{i}}\oplus\fr{m}_{-\xi_{i}})^{\tau}} = \sum_{j_1,\dots,j_{r}}x_{j_1,\dots,j_{r}}B|_{\fr{m}(j_1,\dots,j_r)}
$$
for positive real numbers $x_{\xi}, x_{\xi_i}, x_{j_1,\dots,j_r}$.  The $G$-invariant Riemannian metrics on $M=G/K$ are parametrized by $s$ real positive parameters.
\end{enumerate}
\end{prop}

\subsection{Generalized flag manifolds with $\G_2$-type $\fr{t}$-roots}

A system of positive roots of the Lie group $\G_2$ is given by $\{\alpha _1, \alpha _2, \alpha _1+\alpha _2, 2\alpha _1+\alpha _2, 3\alpha _1+\alpha _2, 3\alpha _1+2\alpha _2\}$,
with heighest root $\widetilde{\alpha}=3\alpha _1+2\alpha _2$.  The corresponding painted Dynkin diagram of the full flag manifold $G_2/T$ is
 
\begin{center}
 \quad 
\small{  \begin{picture}(40,0) (102, 5.9)
\put(106.5, 3.4){\scriptsize $<$}
\put(106,5){\circle*{4.4}}
\put(106, 11.3){\makebox(0,0){$\alpha_1$}}
\put(110, 6.4){\line(1,0){13.8}}
\put(107, 5){\line(1,0){16.5}}
\put(110, 3.9){\line(1,0){14}}
\put(125,5){\circle*{4}}
\put(125.5, 11.3){\makebox(0,0){$\alpha_2$}}
\end{picture} }
\end{center}

From the paper \cite{ArChSa2} we have that  flag manifolds with $\G_2$-type $\fr{t}$-roots system satisfy
$
\Pi\setminus \Pi_{K}=\{\alpha_{i}, \alpha_{j} : \mbox{Mrk}(\alpha_{i})=3, \ \mbox{Mrk}(\alpha_{j})=2\},
$
and are the following:

The highest root $\tilde{\alpha}$ of $F_4$ is given by $\widetilde{\alpha}= 2\alpha_1+4\alpha_2+3\alpha_3+2\alpha_4$, and we have the following painted Dynkin diagram 
  \begin{center}
 \small{  
 \begin{picture}(120,25)(-10,-5)
\put(15, 18) {$\F_4(\al_3, \al_4)$}
\put(16, 0){\circle{4}}
\put(16,8.5){\makebox(0,0){$\alpha_{1}$}}
\put(16,-8){\makebox(0,0){2}}
\put(18, 0){\line(1,0){14}}
\put(32, -2.2){ $<$}
\put(34, 0){\circle{4}}
\put(34, 8.5){\makebox(0,0){$\alpha_{2}$}}
\put(34, -8){\makebox(0,0){4}}
\put(38.5, -1.0){\line(1,0){12}}
\put(38.0, 1.4){\line(1,0){12}}
\put(52, 0){\circle*{4}}
\put(52,8.5){\makebox(0,0){$\alpha_{3}$}}
\put(52,-8){\makebox(0,0){3}}
\put(54, 0){\line(1,0){14}}
\put(70, 0){\circle*{4}}
\put(73,8.5){\makebox(0,0){$\alpha_4$}}
\put(72,-8){\makebox(0,0){2}}
\end{picture}
}
\end{center} 

The highest root $\tilde{\alpha}$ of $\E_6$ is given by $\widetilde{\alpha}=\alpha_1+2\alpha_2+3\alpha_3+2\alpha_4+\alpha_5+2\alpha_6$ and we have the following painted  Dynkin diagram

\begin{center}
\small{ 
\begin{picture}(100,45)(-10,-25)
\put(0, 18) {$\E_6(\al_3, \al_6)$}
\put(0, 0){\circle{4}}
\put(0,8.5){\makebox(0,0){$\alpha_1$}}
\put(0,-8){\makebox(0,0){1}}
\put(2, 0){\line(1,0){14}}
\put(18, 0){\circle{4}}
\put(18,8.5){\makebox(0,0){$\alpha_2$}}
\put(18,-8){\makebox(0,0){2}}
\put(20, 0){\line(1,0){14}}
\put(36, 0){\circle*{4}}
\put(36,8.5){\makebox(0,0){$\alpha_3$}}
\put(32,-8){\makebox(0,0){3}}
\put(38, 0){\line(1,0){14}}
\put(54, 0){\circle{4}}
\put(54,8.5){\makebox(0,0){$\alpha_4$}}
\put(54,-8){\makebox(0,0){2}}
\put(56, 0){\line(1,0){14}}
\put(72, 0){\circle{4}}
\put(72,8.5){\makebox(0,0){$\alpha_5$}}
\put(72,-8){\makebox(0,0){1}}
\put(36, -2){\line(0,-1){14}}
\put(36, -18){\circle*{4}}
\put(43,-25){\makebox(0,0){$\alpha_6$}}
\put(33,-25){\makebox(0,0){2}}
\end{picture} }
\end{center}

The highest root $\tilde{\alpha}$ of $E_7$ is given by $\widetilde{\alpha}=\alpha_1+2\alpha_2+3\alpha_3+4\alpha_4+3\alpha_5+2\alpha_6+2\alpha_7$ and we have the following painted Dynkin diagram

\begin{center}
\small{
\begin{picture}(80, 45)(-10,-22)
\put(-13, 18) {$\E_7(\al_5, \al_6)$}
\put(-18, 0){\circle{4}}
\put(-18,8.5){\makebox(0,0){$\alpha_1$}}
\put(-18,-8){\makebox(0,0){1}}
\put(-16, 0){\line(1,0){14}}
\put(0, 0){\circle{4}}
\put(0,8.5){\makebox(0,0){$\alpha_2$}}
\put(0,-8){\makebox(0,0){2}}
\put(2, 0){\line(1,0){14}}
\put(18, 0){\circle{4}}
\put(18,8.5){\makebox(0,0){$\alpha_3$}}
\put(18,-8){\makebox(0,0){3}}
\put(20, 0){\line(1,0){14}}
\put(36, 0){\circle{4}}
\put(36,8.5){\makebox(0,0){$\alpha_4$}}
\put(32,-8){\makebox(0,0){4}}
\put(38, 0){\line(1,0){14}}
\put(54, 0){\circle*{4}}
\put(54,8.5){\makebox(0,0){$\alpha_5$}}
\put(54,-8){\makebox(0,0){3}}
\put(56, 0){\line(1,0){14}}
\put(72, 0){\circle*{4}}
\put(72,8.5){\makebox(0,0){$\alpha_6$}}
\put(72,-8){\makebox(0,0){2}}
\put(36, -2){\line(0,-1){14}}
\put(36, -18){\circle{4}}
\put(43,-25){\makebox(0,0){$\alpha_7$}}
\put(33,-25){\makebox(0,0){2}}
\end{picture}
 }
\end{center}

The highest root $\tilde{\alpha}$ of $E_8$ is given by $\widetilde{\alpha}= 2\alpha_1+3\alpha_2+4\alpha_3+5\alpha_4+6\alpha_5+4\alpha_6+2\alpha_7+3\alpha_8$ and we have the following painted Dynkin diagram

\begin{center} 
\ \ \ \ \  \small{
\begin{picture}(120,48)(-28.5,-25)
\put(-31, 18) {$\E_8(\al_1, \al_2)$}
\put(-36, 0){\circle*{4}}
\put(-36,8.5){\makebox(0,0){$\alpha_1$}}
\put(-36,-8){\makebox(0,0){2}}
\put(-34, 0){\line(1,0){14}}
\put(-18, 0){\circle*{4}}
\put(-18,8.5){\makebox(0,0){$\alpha_2$}}
\put(-18,-8){\makebox(0,0){3}}
\put(-16, 0){\line(1,0){14}}
\put(0, 0){\circle{4}}
\put(0,8.5){\makebox(0,0){$\alpha_3$}}
\put(0,-8){\makebox(0,0){4}}
\put(2, 0){\line(1,0){14}}
\put(18, 0){\circle{4}}
\put(18,8.5){\makebox(0,0){$\alpha_4$}}
\put(18,-8){\makebox(0,0){5}}
\put(20, 0){\line(1,0){14}}
\put(36, 0){\circle{4}}
\put(36,8.5){\makebox(0,0){$\alpha_5$}}
\put(32,-8){\makebox(0,0){6}}
\put(38, 0){\line(1,0){14}}
\put(54, 0){\circle{4}}
\put(54,8.5){\makebox(0,0){$\alpha_6$}}
\put(54,-8){\makebox(0,0){4}}
\put(56, 0){\line(1,0){14}}
\put(72, 0){\circle{4}}
\put(72,8.5){\makebox(0,0){$\alpha_7$}}
\put(72,-8){\makebox(0,0){2}}
\put(36, -2){\line(0,-1){14}}
\put(36, -18){\circle{4}}
\put(43,-25){\makebox(0,0){$\alpha_8$}}
\put(33,-25){\makebox(0,0){3}}
\end{picture} 
}
\end{center}

For $\Pi \setminus \Pi_{K} =\{\alpha_{i}, \alpha_{j}\}$ we put $\overline{\al}_i = \kappa(\alpha_{i}) $ and $\overline{\al}_j = \kappa(\alpha_{j})$.   We list the sets of all positive $\frak t$-roots $R^{+}_{\frak t}$ in Table 1, which we separate into Type I and Type II.
\\

\begin{center}
{\bf Table 1.} {  Positive $\frak t$-roots $R^{+}_{\frak t}$ for pairs $(\Pi, \Pi_{K})$}  
\end{center}
\begin{center}
\small{
 \begin{tabular}{c|l}
 \hline\hline
  Type I   & Set of all  positive $\frak t$-roots $R^{+}_{\frak t}$    \\
             \hline\hline
    $\F_4(\al_3, \al_4)$   & $ \{  \overline{\al}_3, \ \overline{\al}_4, \  \overline{\al}_3+ \overline{\al}_4,   \ 2\overline{\al}_3+  \overline{\al}_4,  \  3\overline{\al}_3+  \overline{\al}_4,  \ 3 \overline{\al}_3+ 2 \overline{\al}_4 \} $ 
   \\
      \hline
$\E_6(\al_3, \al_6)$  & $ \{  \overline{\al}_3, \ \overline{\al}_6, \  \overline{\al}_3+ \overline{\al}_6,   \ 2 \overline{\al}_3+  \overline{\al}_6,  \  3\overline{\al}_3+ \overline{\al}_6,  \ 3 \overline{\al}_3+  2\overline{\al}_6 \} $  
  \\
       \hline
 $\E_7(\al_5, \al_6)$  &  $\{  \overline{\al}_5, \ \overline{\al}_6, \  \overline{\al}_5+ \overline{\al}_6,   \ 2\overline{\al}_5+  \overline{\al}_6,  \ 3 \overline{\al}_5+  \overline{\al}_6,  \ 3 \overline{\al}_5+ 2 \overline{\al}_6 \} $ \\
       \hline
$\G_2$   &  $  \{  \overline{\al}_1, \ \overline{\al}_2, \  \overline{\al}_1+ \overline{\al}_2,   \ 2\overline{\al}_1+  \overline{\al}_2,  \  3\overline{\al}_1+  \overline{\al}_2,  \ 3 \overline{\al}_1+ 2 \overline{\al}_2 \} $   \\
       \hline \hline      
 Type II & Set of all  positive $\frak t$-roots $R^{+}_{\frak t}$ \\
 \hline\hline   
$\E_8(\al_1, \al_2)$   &  $  \{  \overline{\al}_1, \ \overline{\al}_2, \  \overline{\al}_1+ \overline{\al}_2,   \ \overline{\al}_1+ 2 \overline{\al}_2,  \  \overline{\al}_1+  3 \overline{\al}_2,  \ 2 \overline{\al}_1+ 3 \overline{\al}_2 \} $ \\
       \hline
 \end{tabular}
 }
 \end{center}

From Proposition \ref{isotropy} it is easy to see that the isotropy representation of the above homogeneous spaces is written as a direct sum of six non equivalent $\Ad(K)$-invariant isotropy summands.  For flag manifolds of Type I we set
\begin{equation}\label{compon1}
\begin{array}{lll} 
{\frak m}( 1, 0 ) =  ({\frak m}_{\overline{\al}_{i}} + {\frak m}_{-\overline{\al}_{i}} )^{\tau}, &   &  {\frak m}( 0, 1 ) = ({\frak m}_{\overline{\al}_{j}} + {\frak m}_{-\overline{\al}_{j}})^{\tau}, \\
 {\frak m}( 1, 1 ) =  ({\frak m}_{\overline{\al}_{i}+ \overline{\al}_{j}} + {\frak m}_{-\overline{\al}_{i}- \overline{\al}_{j}})^{\tau}, & & {\frak m}( 2, 1 ) = ({\frak m}_{2\overline{\al}_{i}+ \overline{\al}_{j}} + {\frak m}_{-2\overline{\al}_{i}- \overline{\al}_{j}})^{\tau},\\
  {\frak m}( 3, 1 ) = ({\frak m}_{3\overline{\al}_{i}+ \overline{\al}_{j}} + {\frak m}_{-3\overline{\al}_{i}- \overline{\al}_{j}} )^{\tau},& &  {\frak m}( 3, 2 ) = ({\frak m}_{3\overline{\al}_{i}+ 2\overline{\al}_{j}} + {\frak m}_{-3\overline{\al}_{i}- 2\overline{\al}_{j}})^{\tau},
  \end{array}
  \end{equation} 
and for Type II we set
\begin{equation}\label{compon2}
\begin{array}{lll} 
{\frak m}( 1, 0 ) =  ({\frak m}_{\overline{\al}_{i}} + {\frak m}_{-\overline{\al}_{i}} )^{\tau}, &   &  {\frak m}( 0, 1 ) = ({\frak m}_{\overline{\al}_{j}} + {\frak m}_{-\overline{\al}_{j}})^{\tau}, \\
 {\frak m}( 1, 1 ) =  ({\frak m}_{\overline{\al}_{i}+ \overline{\al}_{j}} + {\frak m}_{-\overline{\al}_{i}- \overline{\al}_{j}})^{\tau}, & & {\frak m}( 1, 2 ) = ({\frak m}_{\overline{\al}_{i}+ 2\overline{\al}_{j}} + {\frak m}_{-\overline{\al}_{i}- 2\overline{\al}_{j}})^{\tau},\\
  {\frak m}( 1, 3 ) = ({\frak m}_{\overline{\al}_{i}+ 3\overline{\al}_{j}} + {\frak m}_{-\overline{\al}_{i}- 3\overline{\al}_{j}} )^{\tau},& &  {\frak m}( 2, 3 ) = ({\frak m}_{2\overline{\al}_{i}+ 3\overline{\al}_{j}} + {\frak m}_{-2\overline{\al}_{i}- 3\overline{\al}_{j}})^{\tau}. 
  \end{array}
  \end{equation} 
By using tables of positive roots (eg. Table B in Appendix of [\cite{FdV}, pp. 528--531]),  we obtain the  dimensions of these spaces as shown in  Table 2. 
\begin{center}
{\bf Table 2.} { Dimensions of irreducible summands with  $G_2$-type $\frak t$-roots}  
\end{center}
 \begin{center}
\small{
 \begin{tabular}{c|c|c|c|c|c|c}
  \hline \hline
  Type I   & $ {\frak m}( 1, 0 )$ & $ {\frak m}( 0, 1 )$ & ${\frak m}( 1, 1 )$ & $ {\frak m}( 2, 1 )$ & $  {\frak m}( 3, 1 )$ & $ {\frak m}( 3, 2 )$  \\
     \hline\hline 
 $\F_4(\al_3, \al_4)$    & 12 & 2 & 12 & 12 &  2 & 2 \\
       \hline 
 $\E_6(\al_3, \al_6)$    & 18 & 2 & 18 & 18 &  2 & 2 \\
       \hline 
$\E_7(\al_5, \al_6)$   & 30 &  2 & 30 & 30 &  2 & 2 \\
       \hline 
 $G_2 $ & 2 & 2 & 2 & 2 &  2 & 2 \\
 \hline\hline
Type II & $ {\frak m}( 1, 0 )$ & $ {\frak m}( 0, 1 )$ & ${\frak m}( 1, 1 )$ & $ {\frak m}( 1, 2 )$ & $  {\frak m}( 1, 3 )$ & $ {\frak m}( 2, 3 )$  \\
\hline\hline
$\E_8(\al_1, \al_2)$    &2 & 54 & 54 & 54  &  2 & 2\\
       \hline 
 \end{tabular}
}
 \end{center}

We consider the generalized flag manifold $M=G/K$ with $\G_2$-type $\fr{t}$-roots.  As we have seen we have the decomposition of $\fr{m}\cong T_{o}(G/K)$ into six irreducible non equivalent $\Ad(K)$-modules as follows:
\begin{eqnarray}
{\rm Type}\ {\rm I} &:& \fr{m} = {\frak m}( 1, 0 ) \oplus {\frak m}( 0, 1 )\oplus{\frak m}( 1, 1 )\oplus {\frak m}( 2, 1 )\oplus {\frak m}( 3, 1 )\oplus {\frak m}( 3, 2 ) \label{dec1}\\
{\rm Type}\ {\rm II} &:& \fr{m} = {\frak m}( 1, 0 )\oplus {\frak m}( 0, 1 )\oplus{\frak m}( 1, 1 )\oplus{\frak m}( 1, 2 )\oplus  {\frak m}( 1, 3 )\oplus{\frak m}( 2, 3 ). \label{dec2}
\end{eqnarray}
For Type I we set $\fr{m}_1 = \fr{m}(1,0), \fr{m}_2=\fr{m}(0,1)$, $\fr{m}_3=\fr{m}(1,1)$, $\fr{m}_4 = \fr{m}(2,1)$, $\fr{m}_5 = \fr{m}(3,1)$, $\fr{m}_6 = \fr{m}(3,2)$, and for Type II we set $\fr{n}_1 = \fr{m}(1,0)$, $\fr{n}_2 = \fr{m}(0,1)$, $\fr{n}_3 = \fr{m}(1,1)$, $\fr{n}_4 = \fr{m}(1,2)$, $\fr{n}_5 = \fr{m}(1,3)$ and $\fr{n}_6 = \fr{m}(2,3)$.  

We now compute the Lie brackets $[\fr{m}_i, \fr{m}_j]$ and $[\fr{n}_i, \fr{n}_j]$ among the real irreducible submodules $\fr{m}_i$ and $\fr{n}_i$ of $\fr{m}$.  According to (\ref{sxesi}), each real submodule $\fr{m}_i$  (or $\fr{n}_i$) associated to the positive $\fr{t}$-root $\xi_i$ can be expressed in terms of root vectors $E_{\pm\al}$ $(\al\in R^{+}_{M})$, such that $\kappa(\al) = \xi_i$.  So from (\ref{LieBracket}) we can compute the brackets $[\fr{m}_i, \fr{m}_j]$ (or $[\fr{n}_i, \fr{n}_j]$), for suitable root vectors $E_\al$.

\begin{lemma}
Let $M=G/K$ be the flag manifold of Type I. Then we obtain that $[\fr{m}_i, \fr{m}_i]\subset\fr{k}$ for $1\leq i\leq 6$, and
\begin{equation}\label{brackets1}
\begin{array}{lllll}
[\fr{m}_1, \fr{m}_2]\subset\fr{m}_3 & [\fr{m}_2, \fr{m}_3]\subset\fr{m}_1 & [\fr{m}_3, \fr{m}_4]\subset\fr{m}_{1}\oplus\fr{m}_6 & [\fr{m}_4, \fr{m}_5]\subset\fr{m}_1\\

[\fr{m}_1, \fr{m}_3]\subset\fr{m}_2\oplus\fr{m}_4 & [\fr{m}_2, \fr{m}_4]\subset\fr{k} & [\fr{m}_3, \fr{m}_5]\subset\fr{k} & [\fr{m}_4, \fr{m}_6]\subset\fr{m}_3\\

[\fr{m}_1, \fr{m}_4]\subset\fr{m}_3\oplus\fr{m}_5 & [\fr{m}_2, \fr{m}_5]\subset\fr{m}_6 & [\fr{m}_3, \fr{m}_6]\subset\fr{m}_4 & [\fr{m}_5, \fr{m}_6]\subset\fr{m}_2\\

[\fr{m}_1, \fr{m}_5]\subset\fr{m}_4 & [\fr{m}_2, \fr{m}_6]\subset\fr{m}_5 & &   \\

[\fr{m}_1, \fr{m}_6]\subset\fr{k}. & & &  
\end{array}
\end{equation}
\end{lemma}

\begin{lemma}
Let $M = G/K$ be a flag manifold of Type II. Then we obtain that $[\fr{n}_i, \fr{n}_i]\subset\fr{k}$ for $1\leq i\leq 6$, and
\begin{equation}\label{brackets2}
\begin{array}{lllll}
[\fr{n}_1, \fr{n}_2]\subset\fr{n}_3 & [\fr{n}_2, \fr{n}_3]\subset\fr{n}_1\oplus\fr{n}_4 & [\fr{n}_3, \fr{n}_4]\subset\fr{n}_{2}\oplus\fr{n}_6 & [\fr{n}_4, \fr{n}_5]\subset\fr{n}_2\\

[\fr{n}_1, \fr{n}_3]\subset\fr{n}_2 & [\fr{n}_2, \fr{n}_4]\subset\fr{n}_3\oplus\fr{n}_5 & [\fr{n}_3, \fr{n}_5]\subset\fr{k} & [\fr{n}_4, \fr{n}_6]\subset\fr{n}_3\\

[\fr{n}_1, \fr{n}_4]\subset\fr{k} & [\fr{n}_2, \fr{n}_5]\subset\fr{n}_4 & [\fr{n}_3, \fr{n}_6]\subset\fr{n}_4 & [\fr{n}_5, \fr{n}_6]\subset\fr{n}_1\\

[\fr{n}_1, \fr{n}_5]\subset\fr{n}_6 & [\fr{n}_2, \fr{n}_6]\subset\fr{k} & &   \\

[\fr{n}_1, \fr{n}_6]\subset\fr{n}_5. & & &  
\end{array}
\end{equation}
\end{lemma}

\section{Equigeodesics}
Let $G/K$ be a generalized flag manifold equipped with a $G$-invariant metric $g$.  It is known that such metrics are in one-to-one correspondence with $\Ad(K)$-invariant scalar products $\langle\cdot,\cdot\rangle$ on $\fr{m}\cong T_{o}(G/K)$ (\cite[Proposition 3.1]{KoNo}).  These in turn correspond $\Ad(K)$-equivariant, positive definite, symmetric operators $\Lambda : \fr{m}\to \fr{m}$ determined by $\langle\cdot, \cdot\rangle = Q(\Lambda\cdot, \cdot)$, where $Q = -B$, the negative of the Killing form on $\fr{g}$.  A curve of the form $\gamma(t) = (\exp tX)\cdot o$ is called {\it equigeodesic} on $G/K$ if it is a geodesic with respect to each invariant metric on $G/K$.  The vector $X$ is called {\it equigeodesic vector}.  The following proposition gives us an algebraic characterization of equigeodesic vectors.

\begin{prop}\textnormal{(\cite{CGN})}
Let $G/K$ be a reductive homogeneous space with reductive decomposition $\fr{g}=\fr{k}\oplus\fr{m}$ and $X\in\fr{m}$ be a non-zero vector.  Then $X$ is a equigeodesic vector if and only if
\begin{equation}\label{equigeod}
[X, \Lambda X]_{\fr{m}}=0,
\end{equation}
for each invariant metric $\Lambda$.
\end{prop}     

To solve equation (\ref{equigeod}) is equivalent to solve a non linear algebraic system of equations whose variables are the coefficients of the vector $X$.  We consider the decomposition $\fr{m} = \sum_{\al\in R_{M}^{+}}\fr{m}_{\al}$ and the basis $\{A_{\al}, B_{\al} : \al\in R_{M}^{+}\}$.  Then, by analysing the Lie brackets  $[A_{\al}, B_{\beta}]$, $[A_{\al}, A_{\beta}]$, $[B_{\al}, B_{\beta}]$ described in (\ref{LieBracket}) it is clear that if the structural constants $N_{\al, \beta}, N_{-\al, \beta}$, $N_{\al, -\beta}$ vanish (e.g. if $\al\pm\beta$ is not a root), then these brackets also vanish and the system can be simplified.  In some cases (depending just on the $\fr{m}_i$-parts of $X$) the nonlinear system vanishes completely (i.e. the system is identically zero).  This motivates the following definition:

\begin{definition}
An equigeodesic vector $X$ is said to be

\smallskip
$(a)$ structural if the algebraic system associated to equation (\ref{equigeod}) vanishes completely.

\smallskip
$(b)$ algebraic if the coordinates of the vector $X$ come from a solution of a (not identically zero) nonlinear algebraic system associated to equation (\ref{equigeod}). 
\end{definition}    

\begin{remark}
\textnormal{From the invariance of the metric $\Lambda$, we have that $\Lambda|_{\fr{m}_i} = \lambda_{i}{\rm Id}_{\fr{m}_i}$, for some $\lambda_i > 0$, for each irreducible component of the isotropy representation.  Therefore, if $X\in\fr{m}_i$ then equation (\ref{equigeod}) is satisfied trivially.}
\end{remark}

We call an equigeodesic vector $X\in\fr{m}$ {\it trivial} if $X\in\fr{m}_i$ for some $i$, otherwise is said to be {\it non trivial}.  It is obvious that the trivial equigeodesic vectors are structural equigeodesic vectors.

\begin{lemma}
Let $G/K$ be a generalized flag manifold with $\G_2$-type $\fr{t}$-roots.  A vector $X=\sum_{i=1}^{6}X_{\fr{m}_i}\in\fr{m} = \bigoplus_{i=1}^{6}\fr{m}_i$ (resp. $X=\sum_{i=1}^{6}X_{\fr{n}_i}\in\fr{n} = \bigoplus_{i=1}^{6}\fr{n}_i$) is equigeodesic if and only if
\begin{equation}\label{lemma}
[X_{\fr{m}_i}, X_{\fr{m}_j}] = 0, \quad (\mbox{resp.}\ [X_{\fr{n}_i}, X_{\fr{n}_j}] = 0)
\end{equation}
where $1\leq i< j\leq 6$.
\end{lemma}
\begin{proof}
If $\pi : \fr{g} \to \fr{m}$ is the projection onto $\fr{m}$, then $\pi([X, \Lambda X])=[X, \Lambda X]_{\fr{m}}$.  Assume that $G/K$ is a flag manifold of Type I.  Let $X=\sum_{i=1}^{6}X_{\fr{m}_i}\in\fr{m} = \bigoplus_{i=1}^{6}\fr{m}_i$.  Then 
\begin{eqnarray*}
&&[X, \Lambda X]_{\fr{m}} = \pi([X, \Lambda X]) = \pi(\big[\sum_{i=1}^{6}X_{\fr{m}_i}, \Lambda(\sum_{i=1}^{6}X_{\fr{m}_{i}})\big]) = \pi(\big[\sum_{i=1}^{6}X_{\fr{m}_{i}}, \sum_{i=1}^{6}\lambda_{i}X_{\fr{m}_{i}}\big])\\
&& = \sum_{i=2}^{6}(\lambda_{i}-\lambda_1)\pi([X_{\fr{m}_1}, X_{\fr{m}_i}]) + \sum_{i=3}^{6}(\lambda_i-\lambda_2)\pi([X_{\fr{m}_2}, X_{\fr{m}_i}]) + \cdots + (\lambda_6 - \lambda_{5})\pi([X_{\fr{m}_{5}}, X_{\fr{m}_6}]) \\
&&= (\lambda_2-\lambda_1)[X_{\fr{m}_1}, X_{\fr{m}_2}] + (\lambda_3-\lambda_1)[X_{\fr{m}_1}, X_{\fr{m}_3}] + (\lambda_4-\lambda_1)[X_{\fr{m}_1}, X_{\fr{m}_4}] + (\lambda_5-\lambda_1)[X_{\fr{m}_1}, X_{\fr{m}_5}] \\
&& + (\lambda_3-\lambda_2)[X_{\fr{m}_2}, X_{\fr{m}_3}] + (\lambda_5-\lambda_2)[X_{\fr{m}_2}, X_{\fr{m}_5}] + (\lambda_6-\lambda_2)[X_{\fr{m}_2}, X_{\fr{m}_6}] + (\lambda_4-\lambda_3)[X_{\fr{m}_3}, X_{\fr{m}_4}]\\
&& + (\lambda_6-\lambda_3)[X_{\fr{m}_3}, X_{\fr{m}_6}] + (\lambda_5-\lambda_4)[X_{\fr{m}_4}, X_{\fr{m}_5}] + (\lambda_6-\lambda_4)[X_{\fr{m}_4}, X_{\fr{m}_6}] + (\lambda_6-\lambda_5)[X_{\fr{m}_5}, X_{\fr{m}_6}].
\end{eqnarray*}
From the bracket relations (\ref{brackets1}) we see that all the above brackets $[X_{\fr{m}_i}, X_{\fr{m}_j}]$ belong to $\fr{m}$.  We know that $X$ is a equigeodesic vector if and only if $[X, \Lambda X]_{\fr{m}} = 0$ for each invariant metric $\Lambda = \{\lambda_1,\ldots,\lambda_6\}$ $(\lambda_1>0,\ldots,\lambda_6>0)$.  This occurs if and only if $[X_{\fr{m}_i}, X_{\fr{m}_j}] = 0$, where $1\leq i<j \leq 6$.  Similarly, if $G/K$ is a flag manifold of Type II we use bracket relations (\ref{brackets2}) and can show that the vector $X=\sum_{i=1}^{6}X_{\fr{n}_i}\in\fr{n} = \bigoplus_{i=1}^{6}\fr{n}_i$ is equigeodesic if and only if $[X_{\fr{n}_i}, X_{\fr{n}_j}] = 0$ where $1\leq i<j \leq 6$. 
\end{proof}

In the papers \cite{GrNe} and \cite{WaZh} the authors gave a family of structural equigeodesic vectors for the generalized flag manifolds $G/K$ with two and $s$ isotropy summands respectively, which depend only on the Lie algebra structure of $\fr{g}$.  
The next theorem provides a family of structural equigeodesic vectors on generalized flag manifolds with $\G_2$-type $\fr{t}$-roots.

\begin{theorem}\label{BasicTheorem}
Let $G/K$ be a generalized flag manifold with $\G_2$-type $\fr{t}$-roots and $\Pi_{K} = \Pi\backslash \{\al_{i_0}, \al_{j_{0}}\}$ of Type I.  Let the positive roots $R^{\fr{m}}(1,0)=\{\beta_{1}^{1}, \ldots, \beta_{k_1}^{1}\}$, $R^{\fr{m}}(0,1)=\{\beta_{1}^{2},\dots, \beta_{k_2}^{2}\}$, $R^{\fr{m}}(1,1)=\{\beta_{1}^{3},\dots, \beta_{k_3}^{3}\}$, $R^{\fr{m}}(2,1)=\{\beta_{1}^{4},\dots, \beta_{k_4}^{4}\}$, $R^{\fr{m}}(3,1)=\{\beta_{1}^{5},\dots, \beta_{k_5}^{5}\}$, $R^{\fr{m}}(3,2)=\{\beta_{1}^{6},\ldots,\beta_{k_6}^{6}\}$.  Suppose that the set $\{\beta_{i_1}^{1},$ $\beta_{i_2}^{2},$ \ldots,$\beta_{i_6}^{6}$ $\ : \ 1\leq i_1\leq k_1,\ldots, 1\leq i_6\leq k_6\}$ satisfies
$$
\beta_{i_p}^{p} \pm \beta_{i_q}^{q} \notin R \ \ \mbox{for}\ \ i_p\neq i_q \ \mbox{and} \ \ 1\leq p< q\leq 6.
$$
Then all vectors in the subspace $\fr{U}_{\beta_{i_1}^{1}}\oplus\fr{U}_{\beta_{i_2}^{2}}\oplus\cdots\oplus\fr{U}_{\beta_{i_6}^{6}}$ are structural equigeodesic vectors.
\end{theorem}
\begin{proof}
Let $X = \sum_{i=1}^{6}X_{\fr{m}_i} = X_{(1,0)}+X_{(0,1)}+X_{(1,1)}+X_{(2,1)}+X_{(3,1)}+X_{(3,2)}$, where 
\begin{eqnarray*}
X_{(1,0)} = \sum_{\al\in R^{\fr{m}}(1,0)}\{\bb{R}A_{\al} + \bb{R}B_{\al}\} = \sum_{\al\in R^{\fr{m}}(1,0)}\fr{U}_{\al}, && X_{(0,1)} = \sum_{\al\in R^{\fr{m}}(0,1)}\{\bb{R}A_{\al} + \bb{R}B_{\al}\} = \sum_{\al\in R^{\fr{m}}(0,1)}\fr{U}_{\al} \\
X_{(1,1)} = \sum_{\al\in R^{\fr{m}}(1,1)}\{\bb{R}A_{\al} + \bb{R}B_{\al}\} = \sum_{\al\in R^{\fr{m}}(1,1)}\fr{U}_{\al},
&& X_{(2,1)} = \sum_{\al\in R^{\fr{m}}(2,1)}\{\bb{R}A_{\al} + \bb{R}B_{\al}\} = \sum_{\al\in R^{\fr{m}}(2,1)}\fr{U}_{\al} \\
X_{(3,1)} = \sum_{\al\in R^{\fr{m}}(3,1)}\{\bb{R}A_{\al} + \bb{R}B_{\al}\} = \sum_{\al\in R^{\fr{m}}(3,1)}\fr{U}_{\al},
&& X_{(3,2)} = \sum_{\al\in R^{\fr{m}}(3,2)}\{\bb{R}A_{\al} + \bb{R}B_{\al}\} = \sum_{\al\in R^{\fr{m}}(3,2)}\fr{U}_{\al}.
\end{eqnarray*}

Since $\beta_{i_{p}}^{p}\pm\beta_{i_{q}}^{q}$ is not a root for $i_p\neq i_q$ and $1\leq p< q\leq 6$, we have that $N_{\beta_{i_{p}}^{p}, \beta_{i_{q}}^{q}} = N_{-\beta_{i_{p}}^{p}, \beta_{i_{q}}^{q}} = N_{\beta_{i_{p}}^{p}, -\beta_{i_{q}}^{q}} = 0$.  A direct computation using the relations of Lemma \ref{basiclemma} shows that the system of equations (\ref{lemma}) vanishes and the vector $X$ is a structural equigeodesic vector.
\end{proof}

Similarly, we can proof the following:

\begin{theorem}\label{BasicTheorem1}
Let $G/K$ be a generalized flag manifold with $\G_2$-type $\fr{t}$-roots and $\Pi_{K} = \Pi\backslash \{\al_{i_0}, \al_{j_{0}}\}$ of Type II.  Let the positive roots $R^{\fr{n}}(1,0)=\{\beta_{1}^{1}, \ldots, \beta_{k_1}^{1}\}$, $R^{\fr{n}}(0,1)=\{\beta_{1}^{2},\dots, \beta_{k_2}^{2}\}$, $R^{\fr{n}}(1,1)=\{\beta_{1}^{3},\dots, \beta_{k_3}^{3}\}$, $R^{\fr{n}}(1,2)=\{\beta_{1}^{4},\dots, \beta_{k_4}^{4}\}$, $R^{\fr{n}}(1,3)=\{\beta_{1}^{5},\dots, \beta_{k_5}^{5}\}$, $R^{\fr{n}}(2,3)=\{\beta_{1}^{6},\ldots,\beta_{k_6}^{6}\}$.  Suppose that the set $\{\beta_{i_1}^{1},$ $\beta_{i_2}^{2},$ \ldots,$\beta_{i_6}^{6}$ $\ : \ 1\leq i_1\leq k_1,\ldots, 1\leq i_6\leq k_6\}$ satisfies
$$
\beta_{i_p}^{p} \pm \beta_{i_q}^{q} \notin R \ \ \mbox{for}\ \ i_p\neq i_q \ \mbox{and} \ \ 1\leq p< q\leq 6.
$$
Then all vectors in the subspace $\fr{U}_{\beta_{i_1}^{1}}\oplus\fr{U}_{\beta_{i_2}^{2}}\oplus\cdots\oplus\fr{U}_{\beta_{i_6}^{6}}$ are structural equigeodesic vectors.
\end{theorem}


\section{Structural Equigeodesic Vectors on Flag Manifolds with $\G_2$-type $\fr{t}$-roots}

In this section we give a family of structural equigeodesic vectors for generalized flag manifolds with $\G_2$-type $\fr{t}$-roots, namely for $\F_4/(\U(3)\times\U(1))$, $\E_6/(\U(3)\times\U(3))$, $\E_7/(\U(6)\times\U(1))$ and $\E_8/(\E_6\times\U(1)\times\U(1))$. We classify the positive roots that satisfy the hypothesis of Theorems \ref{BasicTheorem} and \ref{BasicTheorem1}.  

For the root system of exceptional Lie groups $\F_4, \E_6, \E_7$ and $\E_8$ we use the notation of \cite{AlAr}, where all positive roots are given as linear combinations of the simple roots $\Pi=\{\al_1, \al_2,\ldots,\al_\ell\}$ ($\ell = {\rm rk}\fr{g}^{\bb{C}}$).

\subsection{Structural Equigeodesic vectors on the flag manifold $\F_4/(\U(3)\times\U(1))$} 
Let $\Pi = \{\al_1, \al_2, \al_3, \al_4\}$ be a system of simple roots for $\F_4$ with highest root $\widetilde{\alpha}= 2\alpha_1+4\alpha_2+3\alpha_3+2\alpha_4$. The flag manifold $\F_4/(\U(3)\times\U(1))$ is determined by $\Pi_K = \Pi\backslash \{\al_{3}, \al_{4}\}$.  From Table 1 we have that the the positive $\fr{t}$-roots are given by  $R_{\fr{t}}^{+} = \{  \overline{\al}_3, \ \overline{\al}_4, \  \overline{\al}_3+ \overline{\al}_4,   \ 2\overline{\al}_3+  \overline{\al}_4,  \  3\overline{\al}_3+  \overline{\al}_4,  \ 3 \overline{\al}_3+ 2 \overline{\al}_4 \}$.  According to Proposition \ref{onetone} (1), we obtain the decomposition (\ref{dec1}) where the sumbodules $\fr{m}_i$ are defined by (\ref{compon1}). The sets $R^{\fr{m}}(j_1, j_2) = \{\sum_{i=1}^{4}c_i\al_i\in R_{M}^{+} : c_{3} = j_1, \ c_4 = j_2\}$ are given explicitly as follows:
\begin{eqnarray*}
R^{\fr{m}}(1,0) &=& \{e_3, e_1-e_2, e_3-e_4, e_3+e_4, 1/2(e_1-e_2+e_3+e_4), 1/2(e_1-e_2+e_3-e_4)\} = \{\beta_{1}^1,\ldots, \beta_{6}^{1}\}\\
R^{\fr{m}}(0,1) &=& \{e_2-e_3\} = \{\beta_{1}^{2}\}\\
R^{\fr{m}}(1,1) &=& \{e_2, e_1-e_3, e_2-e_4, e_2+e_4, 1/2(e_1+e_2-e_3+e_4), 1/2(e_1+e_2-e_3-e_4)\} = \{\beta_{1}^{3}, \ldots, \beta_{6}^{3}\}\\
R^{\fr{m}}(2,1) &=& \{e_1, e_1-e_4, e_1+e_4, e_2+e_3, 1/2(e_1+e_2+e_3+e_4), 1/2(e_1+e_2+e_3-e_4)\} = \{\beta_{1}^{4}, \ldots, \beta_{6}^{4}\}\\
R^{\fr{m}}(3,1) &=& \{e_1+e_3\} = \{\beta_{1}^{5}\}\\
R^{\fr{m}}(3,2) &=& \{e_1+e_2\} = \{\beta_{1}^{6}\}. 
\end{eqnarray*}
It is easy to see that the roots which satisfy the hypothesis of Theorem \ref{BasicTheorem} are the following:
\begin{equation*}\label{eq1}
\begin{array}{l}
\beta_{i}^{1}\pm\beta_{j}^{3} \notin R\ \mbox{for every}\ (i, j)\in\{(1,3), (1,4), (2,5), (2,6), (3,1), (3,6), (4,1), (4,5), (5,2), (5,4), (6,2), (6,3)\}\\ \\
\beta_{i}^{1}\pm\beta_{j}^{4} \notin R\ \mbox{for every}\ (i, j)\in\{(1,2), (1,3), (2,5), (2,6), (3,1), (3,5), (4,1), (4,6), (5,2), (5,4), (6,3), (6,4)\}\\ \\
\beta_{i}^{1}\pm\beta_{1}^{6} \notin R \ \mbox{for every}\ i=1,2,\ldots,6\\ \\
\beta_{1}^{2}\pm\beta_{j}^{4} \notin R \ \mbox{for every}\ j=1,2,\ldots,6\\ \\
\beta_{i}^{3}\pm\beta_{j}^{4} \notin R\ \mbox{for every}\ (i,j)\in\{(1,2), (1,3), (2,5), (2,6), (3,1), (3,5), (4,1), (4,6), (5,2), (5,4), (6,3), (6,4)\}\\ \\
\beta_{i}^{3}\pm\beta_{1}^{5} \notin R\ \mbox{for every}\ i=1,2,\ldots,6.
\end{array}
\end{equation*}

From the above roots we can find all the subspaces for which the vectors are structural equigeodesic vectors.  In particular we have the following:
\begin{prop}\label{prop4.1}
The root spaces for the generalized flag manifold $\F_4/(\U(3)\times\U(1))$ whose roots satisfy Theorem \ref{BasicTheorem} are listed in Table 3.  In particular, all vectors in these subspaces are structural equigeodesic vectors. 
\begin{center}
{\bf Table 3.} {{\rm Structural equigeodesic vectors for $\F_4/(\U(3)\times\U(1))$}}  
\end{center} 
\begin{center}
{\small \begin{tabular}{cccc}
\hline\hline
$\fr{U}_{\beta_3^3}\oplus\fr{U}_{\beta_1^1}\oplus\fr{U}_{\beta_{6}^1}$  & 
$\fr{U}_{\beta_1^4}\oplus\fr{U}_{\beta_3^1}\oplus\fr{U}_{\beta_{4}^1}$  &
$\fr{U}_{\beta_3^3}\oplus\fr{U}_{\beta_1^4}\oplus\fr{U}_{\beta_{5}^4}$  &
$\fr{U}_{\beta_4^3}\oplus\fr{U}_{\beta_1^1}\oplus\fr{U}_{\beta_{5}^1}$  \\
$\fr{U}_{\beta_2^4}\oplus\fr{U}_{\beta_1^1}\oplus\fr{U}_{\beta_{5}^1}$ & 
$\fr{U}_{\beta_4^3}\oplus\fr{U}_{\beta_1^4}\oplus\fr{U}_{\beta_{6}^4}$ &
$\fr{U}_{\beta_5^3}\oplus\fr{U}_{\beta_2^1}\oplus\fr{U}_{\beta_{4}^1}$  &
$\fr{U}_{\beta_3^4}\oplus\fr{U}_{\beta_1^1}\oplus\fr{U}_{\beta_{6}^1}$  \\
$\fr{U}_{\beta_5^3}\oplus\fr{U}_{\beta_2^4}\oplus\fr{U}_{\beta_{4}^4}$ &
$\fr{U}_{\beta_6^3}\oplus\fr{U}_{\beta_2^1}\oplus\fr{U}_{\beta_{3}^1}$  &
$\fr{U}_{\beta_4^4}\oplus\fr{U}_{\beta_5^1}\oplus\fr{U}_{\beta_{6}^1}$  &
$\fr{U}_{\beta_6^3}\oplus\fr{U}_{\beta_3^4}\oplus\fr{U}_{\beta_{4}^4}$\\
$\fr{U}_{\beta_1^3}\oplus\fr{U}_{\beta_4^1}\oplus\fr{U}_{\beta_{3}^1}$ &
$\fr{U}_{\beta_5^4}\oplus\fr{U}_{\beta_2^1}\oplus\fr{U}_{\beta_{3}^1}$  & 
$\fr{U}_{\beta_1^4}\oplus\fr{U}_{\beta_3^3}\oplus\fr{U}_{\beta_{4}^3}$ &
$\fr{U}_{\beta_2^3}\oplus\fr{U}_{\beta_5^1}\oplus\fr{U}_{\beta_{6}^1}$  \\
$\fr{U}_{\beta_6^4}\oplus\fr{U}_{\beta_2^1}\oplus\fr{U}_{\beta_{4}^1}$  & 
$\fr{U}_{\beta_2^4}\oplus\fr{U}_{\beta_1^3}\oplus\fr{U}_{\beta_{5}^3}$ &
$\fr{U}_{\beta_1^1}\oplus\fr{U}_{\beta_3^3}\oplus\fr{U}_{\beta_{4}^3}$  &
$\fr{U}_{\beta_1^1}\oplus\fr{U}_{\beta_2^4}\oplus\fr{U}_{\beta_{3}^4}$  \\
$\fr{U}_{\beta_3^4}\oplus\fr{U}_{\beta_1^3}\oplus\fr{U}_{\beta_{6}^3}$ &
$\fr{U}_{\beta_2^1}\oplus\fr{U}_{\beta_5^3}\oplus\fr{U}_{\beta_{6}^3}$ &
$\fr{U}_{\beta_2^1}\oplus\fr{U}_{\beta_5^4}\oplus\fr{U}_{\beta_{6}^4}$ &
$\fr{U}_{\beta_4^4}\oplus\fr{U}_{\beta_5^3}\oplus\fr{U}_{\beta_{6}^3}$ \\
$\fr{U}_{\beta_3^1}\oplus\fr{U}_{\beta_1^3}\oplus\fr{U}_{\beta_{6}^3}$ &
$\fr{U}_{\beta_3^1}\oplus\fr{U}_{\beta_1^4}\oplus\fr{U}_{\beta_{5}^4}$  & 
$\fr{U}_{\beta_5^4}\oplus\fr{U}_{\beta_2^3}\oplus\fr{U}_{\beta_{3}^3}$ &
$\fr{U}_{\beta_4^1}\oplus\fr{U}_{\beta_1^3}\oplus\fr{U}_{\beta_{5}^3}$  \\
$\fr{U}_{\beta_4^1}\oplus\fr{U}_{\beta_1^4}\oplus\fr{U}_{\beta_{6}^4}$ & 
$\fr{U}_{\beta_6^4}\oplus\fr{U}_{\beta_2^3}\oplus\fr{U}_{\beta_{4}^3}$ &
$\fr{U}_{\beta_5^1}\oplus\fr{U}_{\beta_2^3}\oplus\fr{U}_{\beta_{4}^3}$ &
$\fr{U}_{\beta_5^1}\oplus\fr{U}_{\beta_2^4}\oplus\fr{U}_{\beta_{4}^4}$ \\ 
$\fr{U}_{\beta_1^3}\oplus\fr{U}_{\beta_2^4}\oplus\fr{U}_{\beta_{3}^4}$ &
$\fr{U}_{\beta_6^1}\oplus\fr{U}_{\beta_2^3}\oplus\fr{U}_{\beta_{3}^3}$ &
$\fr{U}_{\beta_6^1}\oplus\fr{U}_{\beta_3^4}\oplus\fr{U}_{\beta_{4}^4}$  & 
$\fr{U}_{\beta_2^3}\oplus\fr{U}_{\beta_5^4}\oplus\fr{U}_{\beta_{6}^4}$\\
$\fr{U}_{\beta_1^6}\oplus\fr{U}_{\beta_1^1}\oplus\fr{U}_{\beta_{2}^1}\oplus$ & $\fr{U}_{\beta_1^2}\oplus\fr{U}_{\beta_1^4}\oplus\fr{U}_{\beta_{2}^4}\oplus$  &  $\fr{U}_{\beta_1^5}\oplus\fr{U}_{\beta_1^3}\oplus\fr{U}_{\beta_{2}^3}\oplus$ & 
$\fr{U}_{\beta_{3}^1}\oplus\fr{U}_{\beta_{4}^1}\oplus\fr{U}_{\beta_{5}^1}\oplus\fr{U}_{\beta_{6}^1}$ \\
$\fr{U}_{\beta_3^4}\oplus\fr{U}_{\beta_4^4}\oplus\fr{U}_{\beta_{5}^4}\oplus\fr{U}_{\beta_{6}^4}$ & 
$\fr{U}_{\beta_3^3}\oplus\fr{U}_{\beta_4^3}\oplus\fr{U}_{\beta_{5}^3}\oplus\fr{U}_{\beta_{6}^3}$  & & \\
\hline\hline
\end{tabular} }
\end{center}
\end{prop}

\medskip

\subsection{Structural Equigeodesic vectors on the flag manifold $\E_6/(\U(3)\times\U(3))$} 
Let $\Pi = \{\al_1, \al_2, \al_3, \al_4,$ $\al_5,$ $\al_6\}$ be a system of simple roots for $\E_6$ with highest $\widetilde{\alpha}= \alpha_1+2\alpha_2+3\alpha_3+2\alpha_4+\alpha_5+2\alpha_6$. The flag manifold $\E_6/(\U(3)\times\U(3))$ is determined by $\Pi_K = \Pi\backslash \{\al_{3}, \al_{6}\}$.  From Table 1 we have that the the positive $\fr{t}$-roots are given by  $R_{\fr{t}}^{+} = \{  \overline{\al}_3, \ \overline{\al}_6, \  \overline{\al}_3+ \overline{\al}_6,   \ 2\overline{\al}_3+  \overline{\al}_6,  \  3\overline{\al}_3+  \overline{\al}_6,  \ 3 \overline{\al}_3+ 2 \overline{\al}_6 \}$.  According to Proposition \ref{onetone} (1), we obtain the decomposition (\ref{dec1}) where the sumbodules $\fr{m}_i$ are defined by (\ref{compon1}).  The sets $R^{\fr{m}}(j_1, j_2) = \{\sum_{i=1}^{6}c_i\al_i\in R_{M}^{+} : c_{3} = j_1, \ c_6 = j_2\}$ are given explicitly as follows:
\begin{eqnarray*}
R^{\fr{m}}(1,0) &=& \{e_3-e_4, e_2-e_4, e_3-e_5, e_1-e_4, e_1-e_5, e_1-e_6, e_2-e_5, e_2-e_6, e_3-e_6\} = \{\beta_{1}^1,\ldots, \beta_{9}^{1}\}\\
R^{\fr{m}}(0,1) &=& \{e_4+e_5+e_6+e\} = \{\beta_{1}^{2}\}\\
R^{\fr{m}}(1,1) &=& \{e_1+e_4+e_6+e, e_2+e_4+e_6+e, e_3+e_4+e_5+e, e_3+e_5+e_6+e, e_1+e_4+e_5+e,\\ 
&& e_1+e_5+e_6+e, e_2+e_4+e_5+e, e_2+e_5+e_6+e, e_3+e_4+e_6+e\} = \{\beta_{1}^{3}, \ldots, \beta_{8}^{3}\}\\
R^{\fr{m}}(2,1) &=& \{e_1+e_2+e_5+e, e_1+e_3+e_4+e, e_1+e_3+e_6+e, e_2+e_3+e_4+e, e_2+e_3+e_6+e,\\
&& e_1+e_2+e_4+e, e_1+e_2+e_6+e, e_1+e_3+e_5+e, e_2+e_3+e_5+e\} = \{\beta_{1}^{4}, \ldots, \beta_{8}^{4}\}\\
R^{\fr{m}}(3,1) &=& \{e_1+e_2+e_3+e\} = \{\beta_{1}^{5}\}\\
R^{\fr{m}}(3,2) &=& \{e_1+e_2+e_3+e_4+e_5+e_6+2e\} = \{\beta_{1}^{6}\}. 
\end{eqnarray*}
It is easy to see that all roots which satisfy Theorem \ref{BasicTheorem} are the following:
\begin{equation*}\label{eq2}
\begin{array}{l}
\beta_{i}^{1}\pm\beta_{j}^{3} \notin R\ \mbox{for every}\ (i, j)\in\{(1,3), (1,6), (1,8), (1,9), (2,2), (2,4), (2,6), (2,7), (3,1), (3,2), (3,3), (3,4),\\
(4,1), (4,4), (4,5), (4,8), (5,2), (5,5), (5,6), (5,9), (6,1), (6,3), (6,6), (6,7), (7,1), (7,7), (7,8), (7,9), (8,2),\\
(8,3), (8,5), (8,8), (9,4), (9,5), (9,7), (9,9)\}\\ \\
\beta_{i}^{1}\pm\beta_{j}^{4}\notin R\ \mbox{for every}\ (i, j)\in\{(1,1), (1,2), (1,4), (1,7), (2,3), (2,4), (2,6), (2,8), (3,6), (3,7), (3,8), (3,9),\\
(4,2), (4,5), (4,6), (4,9), (5,1), (5,4), (5,5), (5,8), (6,3), (6,4), (6,7), (6,9), (7,1), (7,2), (7,3), (7,9), (8,2),\\
(8,5), (8,7), (8,8), (9,1), (9,3), (9,5), (9,6)\}\\ \\
\beta_{i}^{1}\pm\beta_{1}^{6}\notin R \ \mbox{for every}\ i=1,2,\ldots,9\\ \\
\beta_{1}^{2}\pm\beta_{j}^{4}\notin R \ \mbox{for every}\ j=1,2,\ldots,9\\ \\
\beta_{i}^{3}\pm\beta_{j}^{4}\notin R \ \mbox{for every}\ (i, j)\in\{(1,1), (1,4), (1,5), (1,8), (1,9), (2,1), (2,2), (2,3), (2,8), (2,9), (3,1), (3,3),\\
(3,5), (3,6), (3,7), (4,1), (4,2), (4,4), (4,6), (4,7), (5,3), (5,4), (5,5), (5,7), (5,9), (6,2), (6,4), (6,5), (6,6),\\
 (6,9), (7,2), (7,3), (7,5), (7,7), (7,8), (8,2), (8,3), (8,4), (8,6), (8,8), (9,1), (9,6), (9,7), (9,8), (9,9)\}\\ \\
\beta_{i}^{3}\pm\beta_{1}^{5}\notin R \ \mbox{for every}\ i=1,2,\ldots,9.
\end{array}
\end{equation*}

From the above roots we can find all the subspaces for which the vectors are structural equigeodesic vectors.  More precisely, from the roots $\beta_{i}^{1}\pm\beta_{j}^{3} \notin R$, $\beta_{i}^{1}\pm\beta_{j}^{4}\notin R$, $\beta_{i}^{1}\pm\beta_{1}^{6}\notin R$, $\beta_{1}^{2}\pm\beta_{j}^{4}\notin R$ and $\beta_{i}^{3}\pm\beta_{1}^{5}\notin R$ we obtain the subspaces in the following tables:
\begin{center}
{\bf Table 4.} { Structural equigeodesic vectors for $\E_6/(\U(3)\times\U(3))$}  
\end{center}  
\begin{center}
{\small \begin{tabular}{cccc}
\hline\hline
$\fr{U}_{\beta_1^3}\oplus\fr{U}_{\beta_3^1}\oplus\fr{U}_{\beta_{4}^1}\oplus\fr{U}_{\beta_6^1}\oplus\fr{U}_{\beta_7^1}$ &
$\fr{U}_{\beta_1^4}\oplus\fr{U}_{\beta_1^1}\oplus\fr{U}_{\beta_{5}^1}\oplus\fr{U}_{\beta_7^1}\oplus\fr{U}_{\beta_9^1}$  &
$\fr{U}_{\beta_1^1}\oplus\fr{U}_{\beta_7^1}\oplus\fr{U}_{\beta_{1}^4}\oplus\fr{U}_{\beta_2^4}$ &
$\fr{U}_{\beta_3^1}\oplus\fr{U}_{\beta_4^1}\oplus\fr{U}_{\beta_{1}^3}\oplus\fr{U}_{\beta_4^3}$
 \\
$\fr{U}_{\beta_2^3}\oplus\fr{U}_{\beta_2^1}\oplus\fr{U}_{\beta_{3}^1}\oplus\fr{U}_{\beta_5^1}\oplus\fr{U}_{\beta_8^1}$ &
$\fr{U}_{\beta_2^4}\oplus\fr{U}_{\beta_1^1}\oplus\fr{U}_{\beta_{4}^1}\oplus\fr{U}_{\beta_7^1}\oplus\fr{U}_{\beta_8^1}$ &
$\fr{U}_{\beta_1^1}\oplus\fr{U}_{\beta_6^1}\oplus\fr{U}_{\beta_{4}^4}\oplus\fr{U}_{\beta_7^4}$ &
$\fr{U}_{\beta_6^1}\oplus\fr{U}_{\beta_7^1}\oplus\fr{U}_{\beta_{1}^3}\oplus\fr{U}_{\beta_7^3}$ \\
$\fr{U}_{\beta_3^3}\oplus\fr{U}_{\beta_1^1}\oplus\fr{U}_{\beta_{3}^1}\oplus\fr{U}_{\beta_6^1}\oplus\fr{U}_{\beta_8^1}$ &
$\fr{U}_{\beta_3^4}\oplus\fr{U}_{\beta_2^1}\oplus\fr{U}_{\beta_{6}^1}\oplus\fr{U}_{\beta_7^1}\oplus\fr{U}_{\beta_9^1}$  &
$\fr{U}_{\beta_3^1}\oplus\fr{U}_{\beta_8^1}\oplus\fr{U}_{\beta_{7}^4}\oplus\fr{U}_{\beta_8^4}$  &
$\fr{U}_{\beta_4^1}\oplus\fr{U}_{\beta_7^1}\oplus\fr{U}_{\beta_{1}^3}\oplus\fr{U}_{\beta_8^3}$\\
$\fr{U}_{\beta_4^3}\oplus\fr{U}_{\beta_2^1}\oplus\fr{U}_{\beta_{3}^1}\oplus\fr{U}_{\beta_4^1}\oplus\fr{U}_{\beta_9^1}$ &
$\fr{U}_{\beta_4^4}\oplus\fr{U}_{\beta_1^1}\oplus\fr{U}_{\beta_{2}^1}\oplus\fr{U}_{\beta_5^1}\oplus\fr{U}_{\beta_6^1}$ &
$\fr{U}_{\beta_4^1}\oplus\fr{U}_{\beta_8^1}\oplus\fr{U}_{\beta_{2}^4}\oplus\fr{U}_{\beta_5^4}$ &
$\fr{U}_{\beta_2^1}\oplus\fr{U}_{\beta_3^1}\oplus\fr{U}_{\beta_{2}^3}\oplus\fr{U}_{\beta_4^3}$ \\
$\fr{U}_{\beta_5^3}\oplus\fr{U}_{\beta_4^1}\oplus\fr{U}_{\beta_{5}^1}\oplus\fr{U}_{\beta_8^1}\oplus\fr{U}_{\beta_9^1}$ &
$\fr{U}_{\beta_5^4}\oplus\fr{U}_{\beta_4^1}\oplus\fr{U}_{\beta_{5}^1}\oplus\fr{U}_{\beta_8^1}\oplus\fr{U}_{\beta_9^1}$  &
$\fr{U}_{\beta_4^1}\oplus\fr{U}_{\beta_9^1}\oplus\fr{U}_{\beta_{5}^4}\oplus\fr{U}_{\beta_6^4}$ &
$\fr{U}_{\beta_5^1}\oplus\fr{U}_{\beta_8^1}\oplus\fr{U}_{\beta_{2}^3}\oplus\fr{U}_{\beta_5^3}$ \\
$\fr{U}_{\beta_6^3}\oplus\fr{U}_{\beta_1^1}\oplus\fr{U}_{\beta_{2}^1}\oplus\fr{U}_{\beta_5^1}\oplus\fr{U}_{\beta_6^1}$ &
$\fr{U}_{\beta_6^4}\oplus\fr{U}_{\beta_2^1}\oplus\fr{U}_{\beta_{3}^1}\oplus\fr{U}_{\beta_4^1}\oplus\fr{U}_{\beta_9^1}$ &
$\fr{U}_{\beta_2^1}\oplus\fr{U}_{\beta_6^1}\oplus\fr{U}_{\beta_{3}^4}\oplus\fr{U}_{\beta_4^4}$ &
$\fr{U}_{\beta_1^1}\oplus\fr{U}_{\beta_7^1}\oplus\fr{U}_{\beta_{8}^3}\oplus\fr{U}_{\beta_9^3}$\\
$\fr{U}_{\beta_7^3}\oplus\fr{U}_{\beta_2^1}\oplus\fr{U}_{\beta_{6}^1}\oplus\fr{U}_{\beta_7^1}\oplus\fr{U}_{\beta_9^1}$ &
$\fr{U}_{\beta_7^4}\oplus\fr{U}_{\beta_1^1}\oplus\fr{U}_{\beta_{3}^1}\oplus\fr{U}_{\beta_6^1}\oplus\fr{U}_{\beta_8^1}$ &
$\fr{U}_{\beta_1^1}\oplus\fr{U}_{\beta_5^1}\oplus\fr{U}_{\beta_{1}^4}\oplus\fr{U}_{\beta_4^4}$  &
$\fr{U}_{\beta_2^1}\oplus\fr{U}_{\beta_5^1}\oplus\fr{U}_{\beta_{2}^3}\oplus\fr{U}_{\beta_6^3}$\\
$\fr{U}_{\beta_8^3}\oplus\fr{U}_{\beta_1^1}\oplus\fr{U}_{\beta_{4}^1}\oplus\fr{U}_{\beta_7^1}\oplus\fr{U}_{\beta_8^1}$ &
$\fr{U}_{\beta_8^4}\oplus\fr{U}_{\beta_2^1}\oplus\fr{U}_{\beta_{3}^1}\oplus\fr{U}_{\beta_5^1}\oplus\fr{U}_{\beta_8^1}$ &
$\fr{U}_{\beta_1^1}\oplus\fr{U}_{\beta_6^1}\oplus\fr{U}_{\beta_{3}^3}\oplus\fr{U}_{\beta_6^3}$ &
$\fr{U}_{\beta_2^1}\oplus\fr{U}_{\beta_6^1}\oplus\fr{U}_{\beta_{2}^3}\oplus\fr{U}_{\beta_7^3}$ \\
$\fr{U}_{\beta_9^3}\oplus\fr{U}_{\beta_1^1}\oplus\fr{U}_{\beta_{5}^1}\oplus\fr{U}_{\beta_7^1}\oplus\fr{U}_{\beta_9^1}$ &
$\fr{U}_{\beta_9^4}\oplus\fr{U}_{\beta_3^1}\oplus\fr{U}_{\beta_{4}^1}\oplus\fr{U}_{\beta_6^1}\oplus\fr{U}_{\beta_7^1}$ &
$\fr{U}_{\beta_1^1}\oplus\fr{U}_{\beta_7^1}\oplus\fr{U}_{\beta_{8}^3}\oplus\fr{U}_{\beta_9^3}$ &
$\fr{U}_{\beta_1^1}\oplus\fr{U}_{\beta_6^1}\oplus\fr{U}_{\beta_{3}^3}\oplus\fr{U}_{\beta_6^3}$ \\
$\fr{U}_{\beta_1^1}\oplus\fr{U}_{\beta_3^3}\oplus\fr{U}_{\beta_{6}^3}\oplus\fr{U}_{\beta_8^3}\oplus\fr{U}_{\beta_9^3}$ &
$\fr{U}_{\beta_1^1}\oplus\fr{U}_{\beta_1^4}\oplus\fr{U}_{\beta_{2}^4}\oplus\fr{U}_{\beta_4^4}\oplus\fr{U}_{\beta_7^4}$ &
$\fr{U}_{\beta_2^1}\oplus\fr{U}_{\beta_6^1}\oplus\fr{U}_{\beta_{6}^3}\oplus\fr{U}_{\beta_7^3}$ &
$\fr{U}_{\beta_1^1}\oplus\fr{U}_{\beta_8^1}\oplus\fr{U}_{\beta_{3}^3}\oplus\fr{U}_{\beta_8^3}$ \\
$\fr{U}_{\beta_2^1}\oplus\fr{U}_{\beta_2^3}\oplus\fr{U}_{\beta_{4}^3}\oplus\fr{U}_{\beta_6^3}\oplus\fr{U}_{\beta_7^3}$ &
$\fr{U}_{\beta_2^1}\oplus\fr{U}_{\beta_3^4}\oplus\fr{U}_{\beta_{4}^4}\oplus\fr{U}_{\beta_6^4}\oplus\fr{U}_{\beta_8^4}$ &
$\fr{U}_{\beta_3^1}\oplus\fr{U}_{\beta_8^1}\oplus\fr{U}_{\beta_{2}^3}\oplus\fr{U}_{\beta_3^3}$ &
$\fr{U}_{\beta_4^1}\oplus\fr{U}_{\beta_9^1}\oplus\fr{U}_{\beta_{4}^3}\oplus\fr{U}_{\beta_5^3}$ \\
$\fr{U}_{\beta_3^1}\oplus\fr{U}_{\beta_1^3}\oplus\fr{U}_{\beta_{2}^3}\oplus\fr{U}_{\beta_3^3}\oplus\fr{U}_{\beta_4^3}$ &
$\fr{U}_{\beta_3^1}\oplus\fr{U}_{\beta_6^4}\oplus\fr{U}_{\beta_{7}^4}\oplus\fr{U}_{\beta_8^4}\oplus\fr{U}_{\beta_9^4}$ &
$\fr{U}_{\beta_4^1}\oplus\fr{U}_{\beta_9^1}\oplus\fr{U}_{\beta_{4}^3}\oplus\fr{U}_{\beta_5^3}$ &
$\fr{U}_{\beta_2^1}\oplus\fr{U}_{\beta_9^1}\oplus\fr{U}_{\beta_{4}^3}\oplus\fr{U}_{\beta_7^3}$ \\
$\fr{U}_{\beta_4^1}\oplus\fr{U}_{\beta_1^3}\oplus\fr{U}_{\beta_{4}^3}\oplus\fr{U}_{\beta_5^3}\oplus\fr{U}_{\beta_8^3}$ &
$\fr{U}_{\beta_4^1}\oplus\fr{U}_{\beta_2^4}\oplus\fr{U}_{\beta_{5}^4}\oplus\fr{U}_{\beta_6^4}\oplus\fr{U}_{\beta_9^4}$ &
$\fr{U}_{\beta_4^1}\oplus\fr{U}_{\beta_8^1}\oplus\fr{U}_{\beta_{5}^3}\oplus\fr{U}_{\beta_8^3}$ &
$\fr{U}_{\beta_4^1}\oplus\fr{U}_{\beta_8^1}\oplus\fr{U}_{\beta_{5}^3}\oplus\fr{U}_{\beta_8^3}$ \\
$\fr{U}_{\beta_5^1}\oplus\fr{U}_{\beta_2^3}\oplus\fr{U}_{\beta_{5}^3}\oplus\fr{U}_{\beta_6^3}\oplus\fr{U}_{\beta_9^3}$ &
$\fr{U}_{\beta_5^1}\oplus\fr{U}_{\beta_1^4}\oplus\fr{U}_{\beta_{4}^4}\oplus\fr{U}_{\beta_5^4}\oplus\fr{U}_{\beta_8^4}$ &
$\fr{U}_{\beta_1^1}\oplus\fr{U}_{\beta_5^1}\oplus\fr{U}_{\beta_{6}^3}\oplus\fr{U}_{\beta_9^3}$ &
$\fr{U}_{\beta_5^1}\oplus\fr{U}_{\beta_9^1}\oplus\fr{U}_{\beta_{5}^3}\oplus\fr{U}_{\beta_9^3}$ \\
$\fr{U}_{\beta_6^1}\oplus\fr{U}_{\beta_1^3}\oplus\fr{U}_{\beta_{3}^3}\oplus\fr{U}_{\beta_6^3}\oplus\fr{U}_{\beta_7^3}$ &
$\fr{U}_{\beta_6^1}\oplus\fr{U}_{\beta_3^4}\oplus\fr{U}_{\beta_{4}^4}\oplus\fr{U}_{\beta_7^4}\oplus\fr{U}_{\beta_9^4}$ & 
$\fr{U}_{\beta_3^1}\oplus\fr{U}_{\beta_6^1}\oplus\fr{U}_{\beta_{1}^3}\oplus\fr{U}_{\beta_3^3}$  &
$\fr{U}_{\beta_2^1}\oplus\fr{U}_{\beta_6^1}\oplus\fr{U}_{\beta_{6}^3}\oplus\fr{U}_{\beta_7^3}$\\
$\fr{U}_{\beta_7^1}\oplus\fr{U}_{\beta_1^3}\oplus\fr{U}_{\beta_{7}^3}\oplus\fr{U}_{\beta_8^3}\oplus\fr{U}_{\beta_9^3}$ &
$\fr{U}_{\beta_7^1}\oplus\fr{U}_{\beta_1^4}\oplus\fr{U}_{\beta_{2}^4}\oplus\fr{U}_{\beta_3^4}\oplus\fr{U}_{\beta_9^4}$ & 
$\fr{U}_{\beta_{1}^{6}}\bigoplus_{i=1}^{9}\fr{U}_{\beta_{i}^{1}}$  &
$\fr{U}_{\beta_7^1}\oplus\fr{U}_{\beta_9^1}\oplus\fr{U}_{\beta_{7}^3}\oplus\fr{U}_{\beta_9^3}$\\
$\fr{U}_{\beta_8^1}\oplus\fr{U}_{\beta_2^3}\oplus\fr{U}_{\beta_{3}^3}\oplus\fr{U}_{\beta_5^3}\oplus\fr{U}_{\beta_8^3}$ &
$\fr{U}_{\beta_8^1}\oplus\fr{U}_{\beta_2^4}\oplus\fr{U}_{\beta_{5}^4}\oplus\fr{U}_{\beta_7^4}\oplus\fr{U}_{\beta_8^4}$ & 
$\fr{U}_{\beta_{1}^{2}}\bigoplus_{i=1}^{9}\fr{U}_{\beta_{i}^{4}}$ & 
$\fr{U}_{\beta_1^1}\oplus\fr{U}_{\beta_7^1}\oplus\fr{U}_{\beta_{1}^4}\oplus\fr{U}_{\beta_2^4}$\\
$\fr{U}_{\beta_9^1}\oplus\fr{U}_{\beta_4^3}\oplus\fr{U}_{\beta_{5}^3}\oplus\fr{U}_{\beta_7^3}\oplus\fr{U}_{\beta_9^3}$ &
$\fr{U}_{\beta_9^1}\oplus\fr{U}_{\beta_1^4}\oplus\fr{U}_{\beta_{3}^4}\oplus\fr{U}_{\beta_5^4}\oplus\fr{U}_{\beta_6^4}$  &
$\fr{U}_{\beta_{1}^{5}}\bigoplus_{i=1}^{9}\fr{U}_{\beta_{i}^{3}}$ & 
$\fr{U}_{\beta_7^1}\oplus\fr{U}_{\beta_9^1}\oplus\fr{U}_{\beta_{1}^4}\oplus\fr{U}_{\beta_3^4}$
\\
%

$\fr{U}_{\beta_1^1}\oplus\fr{U}_{\beta_5^1}\oplus\fr{U}_{\beta_{1}^4}\oplus\fr{U}_{\beta_4^4}$ &
$\fr{U}_{\beta_1^1}\oplus\fr{U}_{\beta_8^1}\oplus\fr{U}_{\beta_{2}^4}\oplus\fr{U}_{\beta_7^4}$ &
$\fr{U}_{\beta_2^1}\oplus\fr{U}_{\beta_9^1}\oplus\fr{U}_{\beta_{3}^4}\oplus\fr{U}_{\beta_6^4}$ &
$\fr{U}_{\beta_2^1}\oplus\fr{U}_{\beta_5^1}\oplus\fr{U}_{\beta_{4}^4}\oplus\fr{U}_{\beta_8^4}$\\
%
$\fr{U}_{\beta_5^1}\oplus\fr{U}_{\beta_9^1}\oplus\fr{U}_{\beta_{1}^4}\oplus\fr{U}_{\beta_5^4}$ &
$\fr{U}_{\beta_4^1}\oplus\fr{U}_{\beta_7^1}\oplus\fr{U}_{\beta_{2}^4}\oplus\fr{U}_{\beta_9^4}$ &
$\fr{U}_{\beta_6^1}\oplus\fr{U}_{\beta_7^1}\oplus\fr{U}_{\beta_{3}^4}\oplus\fr{U}_{\beta_9^4}$ &
$\fr{U}_{\beta_4^1}\oplus\fr{U}_{\beta_9^1}\oplus\fr{U}_{\beta_{5}^4}\oplus\fr{U}_{\beta_6^4}$\\
%
$\fr{U}_{\beta_4^1}\oplus\fr{U}_{\beta_8^1}\oplus\fr{U}_{\beta_{2}^4}\oplus\fr{U}_{\beta_5^4}$ &
$\fr{U}_{\beta_2^1}\oplus\fr{U}_{\beta_6^1}\oplus\fr{U}_{\beta_{3}^4}\oplus\fr{U}_{\beta_4^4}$ &
$\fr{U}_{\beta_1^1}\oplus\fr{U}_{\beta_6^1}\oplus\fr{U}_{\beta_{4}^4}\oplus\fr{U}_{\beta_7^4}$ &
$\fr{U}_{\beta_2^1}\oplus\fr{U}_{\beta_3^1}\oplus\fr{U}_{\beta_{6}^4}\oplus\fr{U}_{\beta_8^4}$\\
$\fr{U}_{\beta_3^1}\oplus\fr{U}_{\beta_8^1}\oplus\fr{U}_{\beta_{7}^4}\oplus\fr{U}_{\beta_8^4}$ &
$\fr{U}_{\beta_3^1}\oplus\fr{U}_{\beta_6^1}\oplus\fr{U}_{\beta_{7}^4}\oplus\fr{U}_{\beta_9^4}$ & &\\
\hline\hline
\end{tabular} }
\end{center}

Now from the roots $\beta_{i}^{3}\pm\beta_{j}^{4}\notin R$ we have the following subspaces:
\begin{center}
{\bf Table 5.} { Structural equigeodesic vectors for $\E_6/(\U(3)\times\U(3))$}  
\end{center}  
\begin{center}
{\small \begin{tabular}{cccc}
\hline\hline
$\fr{U}_{\beta_1^4}\oplus\fr{U}_{\beta_1^3}\oplus\fr{U}_{\beta_{2}^3}\oplus\fr{U}_{\beta_3^3}\oplus\fr{U}_{\beta_4^3}\oplus\fr{U}_{\beta_9^3}$ & 
$\fr{U}_{\beta_2^4}\oplus\fr{U}_{\beta_3^4}\oplus\fr{U}_{\beta_{8}^4}\oplus\fr{U}_{\beta_2^3}\oplus\fr{U}_{\beta_7^3}\oplus\fr{U}_{\beta_8^3}$ & 
$\fr{U}_{\beta_4^4}\oplus\fr{U}_{\beta_7^4}\oplus\fr{U}_{\beta_{4}^3}\oplus\fr{U}_{\beta_5^3}$
\\
$\fr{U}_{\beta_2^4}\oplus\fr{U}_{\beta_2^3}\oplus\fr{U}_{\beta_{4}^3}\oplus\fr{U}_{\beta_6^3}\oplus\fr{U}_{\beta_7^3}\oplus\fr{U}_{\beta_8^3}$ & 
$\fr{U}_{\beta_4^4}\oplus\fr{U}_{\beta_5^4}\oplus\fr{U}_{\beta_{9}^4}\oplus\fr{U}_{\beta_1^3}\oplus\fr{U}_{\beta_5^3}\oplus\fr{U}_{\beta_6^3}$ & 
$\fr{U}_{\beta_4^4}\oplus\fr{U}_{\beta_8^4}\oplus\fr{U}_{\beta_{1}^3}\oplus\fr{U}_{\beta_8^3}$
\\
$\fr{U}_{\beta_3^4}\oplus\fr{U}_{\beta_2^3}\oplus\fr{U}_{\beta_{3}^3}\oplus\fr{U}_{\beta_5^3}\oplus\fr{U}_{\beta_7^3}\oplus\fr{U}_{\beta_8^3}$ & 
$\fr{U}_{\beta_1^4}\oplus\fr{U}_{\beta_6^4}\oplus\fr{U}_{\beta_{7}^4}\oplus\fr{U}_{\beta_3^3}\oplus\fr{U}_{\beta_4^3}\oplus\fr{U}_{\beta_9^3}$ & 
$\fr{U}_{\beta_5^4}\oplus\fr{U}_{\beta_6^4}\oplus\fr{U}_{\beta_{3}^3}\oplus\fr{U}_{\beta_6^3}$
\\
$\fr{U}_{\beta_4^4}\oplus\fr{U}_{\beta_1^3}\oplus\fr{U}_{\beta_{5}^3}\oplus\fr{U}_{\beta_6^3}\oplus\fr{U}_{\beta_8^3}\oplus\fr{U}_{\beta_4^3}$ & 
$\fr{U}_{\beta_2^4}\oplus\fr{U}_{\beta_4^4}\oplus\fr{U}_{\beta_{6}^4}\oplus\fr{U}_{\beta_6^3}\oplus\fr{U}_{\beta_8^3}\oplus\fr{U}_{\beta_4^3}$ & 
$\fr{U}_{\beta_5^4}\oplus\fr{U}_{\beta_8^4}\oplus\fr{U}_{\beta_{1}^3}\oplus\fr{U}_{\beta_7^3}$
\\
$\fr{U}_{\beta_5^4}\oplus\fr{U}_{\beta_1^3}\oplus\fr{U}_{\beta_{3}^3}\oplus\fr{U}_{\beta_5^3}\oplus\fr{U}_{\beta_6^3}\oplus\fr{U}_{\beta_7^3}$ & 
$\fr{U}_{\beta_1^4}\oplus\fr{U}_{\beta_8^4}\oplus\fr{U}_{\beta_{9}^4}\oplus\fr{U}_{\beta_1^3}\oplus\fr{U}_{\beta_2^3}\oplus\fr{U}_{\beta_9^3}$ & 
$\fr{U}_{\beta_6^4}\oplus\fr{U}_{\beta_8^4}\oplus\fr{U}_{\beta_{8}^3}\oplus\fr{U}_{\beta_9^3}$
\\
$\fr{U}_{\beta_6^4}\oplus\fr{U}_{\beta_3^3}\oplus\fr{U}_{\beta_{4}^3}\oplus\fr{U}_{\beta_6^3}\oplus\fr{U}_{\beta_8^3}\oplus\fr{U}_{\beta_9^3}$  & 
$\fr{U}_{\beta_3^4}\oplus\fr{U}_{\beta_5^4}\oplus\fr{U}_{\beta_{7}^4}\oplus\fr{U}_{\beta_3^3}\oplus\fr{U}_{\beta_5^3}\oplus\fr{U}_{\beta_7^3}$ &
$\fr{U}_{\beta_6^4}\oplus\fr{U}_{\beta_9^4}\oplus\fr{U}_{\beta_{6}^3}\oplus\fr{U}_{\beta_9^3}$
\\
$\fr{U}_{\beta_7^4}\oplus\fr{U}_{\beta_3^3}\oplus\fr{U}_{\beta_{4}^3}\oplus\fr{U}_{\beta_5^3}\oplus\fr{U}_{\beta_7^3}\oplus\fr{U}_{\beta_9^3}$  & 
$\fr{U}_{\beta_1^3}\oplus\fr{U}_{\beta_5^3}\oplus\fr{U}_{\beta_{6}^3}\oplus\fr{U}_{\beta_5^4}\oplus\fr{U}_{\beta_9^4}$ &
$\fr{U}_{\beta_7^4}\oplus\fr{U}_{\beta_8^4}\oplus\fr{U}_{\beta_{7}^3}\oplus\fr{U}_{\beta_9^3}$
\\
$\fr{U}_{\beta_8^4}\oplus\fr{U}_{\beta_1^3}\oplus\fr{U}_{\beta_{2}^3}\oplus\fr{U}_{\beta_7^3}\oplus\fr{U}_{\beta_8^3}\oplus\fr{U}_{\beta_9^3}$  & 
$\fr{U}_{\beta_3^3}\oplus\fr{U}_{\beta_4^3}\oplus\fr{U}_{\beta_{9}^3}\oplus\fr{U}_{\beta_6^4}\oplus\fr{U}_{\beta_7^4}$ &
$\fr{U}_{\beta_7^4}\oplus\fr{U}_{\beta_9^4}\oplus\fr{U}_{\beta_{5}^3}\oplus\fr{U}_{\beta_9^3}$
\\
$\fr{U}_{\beta_9^4}\oplus\fr{U}_{\beta_1^3}\oplus\fr{U}_{\beta_{2}^3}\oplus\fr{U}_{\beta_5^3}\oplus\fr{U}_{\beta_6^3}\oplus\fr{U}_{\beta_9^3}$ & 
$\fr{U}_{\beta_1^4}\oplus\fr{U}_{\beta_2^4}\oplus\fr{U}_{\beta_{2}^3}\oplus\fr{U}_{\beta_4^3}$ & $\fr{U}_{\beta_1^4}\oplus\fr{U}_{\beta_3^4}\oplus\fr{U}_{\beta_{2}^3}\oplus\fr{U}_{\beta_3^3}$ \\
$\fr{U}_{\beta_1^3}\oplus\fr{U}_{\beta_1^4}\oplus\fr{U}_{\beta_{4}^4}\oplus\fr{U}_{\beta_5^4}\oplus\fr{U}_{\beta_8^4}\oplus\fr{U}_{\beta_9^4}$ & $\fr{U}_{\beta_2^3}\oplus\fr{U}_{\beta_1^4}\oplus\fr{U}_{\beta_{2}^4}\oplus\fr{U}_{\beta_3^4}\oplus\fr{U}_{\beta_8^4}\oplus\fr{U}_{\beta_9^4}$ & $\fr{U}_{\beta_1^4}\oplus\fr{U}_{\beta_4^4}\oplus\fr{U}_{\beta_{1}^3}\oplus\fr{U}_{\beta_4^3}$
\\
$\fr{U}_{\beta_3^3}\oplus\fr{U}_{\beta_1^4}\oplus\fr{U}_{\beta_{3}^4}\oplus\fr{U}_{\beta_5^4}\oplus\fr{U}_{\beta_6^4}\oplus\fr{U}_{\beta_7^4}$ & 
$\fr{U}_{\beta_1^4}\oplus\fr{U}_{\beta_5^4}\oplus\fr{U}_{\beta_{1}^3}\oplus\fr{U}_{\beta_3^3}$ & $\fr{U}_{\beta_2^4}\oplus\fr{U}_{\beta_5^4}\oplus\fr{U}_{\beta_{6}^3}\oplus\fr{U}_{\beta_7^3}$
\\
$\fr{U}_{\beta_4^3}\oplus\fr{U}_{\beta_1^4}\oplus\fr{U}_{\beta_{2}^4}\oplus\fr{U}_{\beta_4^4}\oplus\fr{U}_{\beta_6^4}\oplus\fr{U}_{\beta_7^4}$  & $\fr{U}_{\beta_5^3}\oplus\fr{U}_{\beta_3^4}\oplus\fr{U}_{\beta_{4}^4}\oplus\fr{U}_{\beta_5^4}\oplus\fr{U}_{\beta_7^4}\oplus\fr{U}_{\beta_9^4}$ & $\fr{U}_{\beta_2^4}\oplus\fr{U}_{\beta_7^4}\oplus\fr{U}_{\beta_{4}^3}\oplus\fr{U}_{\beta_9^3}$ 
\\
$\fr{U}_{\beta_6^3}\oplus\fr{U}_{\beta_2^4}\oplus\fr{U}_{\beta_{4}^4}\oplus\fr{U}_{\beta_5^4}\oplus\fr{U}_{\beta_6^4}\oplus\fr{U}_{\beta_9^4}$ & 
$\fr{U}_{\beta_2^4}\oplus\fr{U}_{\beta_9^4}\oplus\fr{U}_{\beta_{2}^3}\oplus\fr{U}_{\beta_6^3}$  & $\fr{U}_{\beta_3^4}\oplus\fr{U}_{\beta_4^4}\oplus\fr{U}_{\beta_{5}^3}\oplus\fr{U}_{\beta_8^3}$ 
\\
$\fr{U}_{\beta_7^3}\oplus\fr{U}_{\beta_2^4}\oplus\fr{U}_{\beta_{3}^4}\oplus\fr{U}_{\beta_5^4}\oplus\fr{U}_{\beta_7^4}\oplus\fr{U}_{\beta_8^4}$ &  $\fr{U}_{\beta_8^3}\oplus\fr{U}_{\beta_2^4}\oplus\fr{U}_{\beta_{3}^4}\oplus\fr{U}_{\beta_4^4}\oplus\fr{U}_{\beta_6^4}\oplus\fr{U}_{\beta_8^4}$ & $\fr{U}_{\beta_3^4}\oplus\fr{U}_{\beta_6^4}\oplus\fr{U}_{\beta_{3}^3}\oplus\fr{U}_{\beta_8^3}$ 
\\
$\fr{U}_{\beta_9^3}\oplus\fr{U}_{\beta_1^4}\oplus\fr{U}_{\beta_{6}^4}\oplus\fr{U}_{\beta_7^4}\oplus\fr{U}_{\beta_8^4}\oplus\fr{U}_{\beta_9^4}$ & 
$\fr{U}_{\beta_3^4}\oplus\fr{U}_{\beta_9^4}\oplus\fr{U}_{\beta_{2}^3}\oplus\fr{U}_{\beta_5^3}$ &
\\
\hline\hline
\end{tabular} }
\end{center}

Hence we obtain the following:
\begin{prop}\label{prop4.2}
The root spaces for the generalized flag manifold $\E_6(\al_3, \al_6)=\E_6/(\U(3)\times\U(3))$, with all vectors are structural equigeodesic vectors are described in Tables 4 and 5.
\end{prop}

\subsection{Structural Equigeodesic vectors on the flag manifold $\E_7/(\U(6)\times\U(1))$} 
Let $\Pi = \{\al_1, \al_2, \al_3, \al_4,$ $\al_5,$ $\al_6, \al_7\}$ be a system of simple roots for $\E_7$ with highest $\widetilde{\alpha}= \alpha_1+2\alpha_2+3\alpha_3+4\alpha_4+3\alpha_5+2\alpha_6+2\alpha_7$. The flag manifold $\E_7/(\U(6)\times\U(1))$ is determined by $\Pi_K = \Pi\backslash \{\al_{5}, \al_{6}\}$.  From Table 1 we have that the the positive $\fr{t}$-roots are given by  $R_{\fr{t}}^{+} = \{  \overline{\al}_5, \ \overline{\al}_6, \  \overline{\al}_5+ \overline{\al}_6,   \ 2\overline{\al}_5+  \overline{\al}_6,  \  3\overline{\al}_5+  \overline{\al}_6,  \ 3 \overline{\al}_5+ 2 \overline{\al}_6 \}$.  According to Proposition \ref{onetone} (1), we obtain the decomposition (\ref{dec1}) where the sumbodules $\fr{m}_i$ are defined by (\ref{compon1}).  The sets $R^{\fr{m}}(j_1, j_2) = \{\sum_{i=1}^{7}c_i\al_i\in R_{M}^{+} : c_{5} = j_1, \ c_6 = j_2\}$ are given explicitly as follows:
\begin{eqnarray*}
R^{\fr{m}}(1,0) &=& \{e_1-e_6, e_2-e_6, e_3-e_6, e_4-e_6, e_5-e_6, e_4+e_5+e_7+e_8, e_3+e_5+e_7+e_8,\\ 
&& e_3+e_4+e_7+e_8, e_2+e_5+e_7+e_8, e_2+e_2+e_7+e_8, e_2+e_3+e_7+e_8, e_1+e_5+e_7+e_8,\\
&& e_1+e_4+e_7+e_8, e_1+e_3+e_7+e_8, e_1+e_2+e_7+e_8\} = \{\beta_{1}^1,\ldots, \beta_{15}^{1}\}\\
R^{\fr{m}}(0,1) &=& \{e_6-e_7\} = \{\beta_{1}^{2}\}\\
R^{\fr{m}}(1,1) &=& \{e_1-e_7, e_2-e_7, e_3-e_7, e_4-e_7, e_5-e_7, e_4+e_5+e_6+e_8, e_3+e_5+e_6+e_8,\\
&& e_3+e_4+e_6+e_8, e_2+e_5+e_6+e_8, e_2+e_4+e_6+e_8, e_2+e_3+e_6+e_8, e_1+e_5+e_6+e_8,\\
&& e_1+e_4+e_6+e_8, e_1+e_3+e_6+e_8, e_1+e_2+e_6+e_8\} = \{\beta_{1}^{3}, \ldots, \beta_{15}^{3}\}\\
R^{\fr{m}}(2,1) &=& \{e_3+e_4+e_5+e_8, e_2+e_4+e_5+e_8, e_2+e_3+e_5+e_8, e_2+e_3+e_4+e_8,\\
&& e_1+e_4+e_5+e_8, e_1+e_3+e_5+e_8, e_1+e_3+e_4+e_8, e_1+e_2+e_5+e_8, e_1+e_2+e_4+e_8, \\
&& e_1+e_2+e_3+e_8, -(e_1-e_8), -(e_2-e_8), -(e_3-e_8), -(e_4-e_8), -(e_5-e_8)\} = \{\beta_{1}^{4}, \ldots, \beta_{15}^{4}\}\\
R^{\fr{m}}(3,1) &=& \{-(e_6-e_8)\} = \{\beta_{1}^{5}\}\\
R^{\fr{m}}(3,2) &=& \{-(e_7-e_8)\} = \{\beta_{1}^{6}\}. 
\end{eqnarray*}

The roots which satisfy Theorem \ref{BasicTheorem} are the following:\\
$\beta_{1}^{1}\pm\beta_{j}^{3}, j =2,3,4,5,12,13,14,15;$  $\beta_{2}^1\pm\beta_{j}^3,  j =1,3,4,5,9,10,11,15;$  $\beta_{3}^{1}\pm\beta_{j}^{3}, j =1,2,4,5,7,8,11,14;$ 
$\beta_{3}^{1}\pm\beta_{j}^{3}, j =1,2,4,5,7,8,11,14;$ $\beta_{1}^{4}\pm\beta_{j}^{3}, j =1,2,3,5,6,8,10,13;$
$\beta_{5}^{1}\pm\beta_{j}^{3}, j =1,2,3,4,6,7,9,12;$ $\beta_{6}^{1}\pm\beta_{j}^{3}, j =4,5,7,8,9,10,12,13;$
$\beta_{7}^{1}\pm\beta_{j}^{3}, j =3,5,6,8,9,11,12,14;$ $\beta_{8}^{1}\pm\beta_{j}^{3}, j =3,4,6,7,10,11,13,14;$
$\beta_{9}^{1}\pm\beta_{j}^{3}, j =2,5,6,7,10,11,12,15;$ $\beta_{10}^{1}\pm\beta_{j}^{3}, j =2,4,6,8,9,11,13,15;$ 
$\beta_{11}^{1}\pm\beta_{j}^{3}, j =2,3,7,8,9,10,14,15;$ $\beta_{12}^{1}\pm\beta_{j}^{3}, j =1,5,6,7,9,13,14,15;$ 
$\beta_{13}^{1}\pm\beta_{j}^{3}, j =1,4,5,6,8,10,12,14,15;$  $\beta_{14}^{1}\pm\beta_{j}^{3}, j =1,3,7,8,11,12,13,15;$ 
$\beta_{15}^{1}\pm\beta_{j}^{3}, j =1,2,9,10,11,12,13,14$  $\beta_{1}^{1}\pm\beta_{j}^{4}, j =1,2,3,4,12,13,14,15;$
$\beta_{2}^{1}\pm\beta_{j}^{4}, j =1,5,6,7,11,13,14,15;$ $\beta_{3}^{1}\pm\beta_{j}^{4}, j =2,5,8,9,11,12,14,15;$
$\beta_{4}^{1}\pm\beta_{j}^{4}, j =3,6,8,10,11,12,13,15;$ $\beta_{5}^{1}\pm\beta_{j}^{4}, j =4,7,9,10,11,12,13,14;$
$\beta_{6}^{1}\pm\beta_{j}^{4}, j =3,4,6,7,8,9,14,15;$ $\beta_{7}^{1}\pm\beta_{j}^{4}, j =2,4,5,7,8,10,13,15;$
$\beta_{8}^{1}\pm\beta_{j}^{4}, j =2,3,5,6,9,10,13,14;$ $\beta_{9}^{1}\pm\beta_{j}^{4}, j =1,4,5,6,9,10,12,15;$
$\beta_{10}^{1}\pm\beta_{j}^{4}, j =1,3,5,7,8,10,12,4;$ $\beta_{11}^{1}\pm\beta_{j}^{4}, j =1,2,6,7,8,9,12,13;$
$\beta_{12}^{1}\pm\beta_{j}^{4}, j =1,2,3,7,9,10,11,15;$ $\beta_{13}^{1}\pm\beta_{j}^{4}, j =1,2,4,6,8,10,11,14;$
$\beta_{14}^{1}\pm\beta_{j}^{4}, j =1,3,4,5,8,9,11,13;$ $\beta_{15}^{1}\pm\beta_{j}^{4}, j =2,3,4,5,6,7,11,12;$
$\beta_{1}^{3}\pm\beta_{j}^{4}, j =1,2,3,4,12,13,14,15;$ $\beta_{2}^{3}\pm\beta_{j}^{4}, j =1,5,6,7,11,13,14,15;$
$\beta_{3}^{3}\pm\beta_{j}^{4}, j =2,5,8,9,11,12,14,15;$ $\beta_{4}^{3}\pm\beta_{j}^{4}, j =3,6,8,10,11,12,13,15;$
$\beta_{5}^{3}\pm\beta_{j}^{4}, j =4,7,9,10,11,12,13,14;$ $\beta_{6}^{3}\pm\beta_{j}^{4}, j =3,4,6,7,8,9,14,15;$
$\beta_{7}^{3}\pm\beta_{j}^{4}, j =2,4,5,7,8,10,13,15;$ $\beta_{8}^{3}\pm\beta_{j}^{4}, j =2,3,5,6,9,10,13,14;$
$\beta_{9}^{3}\pm\beta_{j}^{4}, j =1,4,5,6,9,10,12,15;$ $\beta_{10}^{3}\pm\beta_{j}^{4}, j =1,3,5,7,8,10,12,14;$
$\beta_{11}^{3}\pm\beta_{j}^{4}, j =1,2,6,7,8,9,12,13;$ $\beta_{12}^{3}\pm\beta_{j}^{4}, j =1,2,3,7,9,10,11,15;$
$\beta_{13}^{3}\pm\beta_{j}^{4}, j =1,2,4,6,8,10,11,14;$ $\beta_{14}^{3}\pm\beta_{j}^{4}, j =1,3,4,5,8,9,11,13;$
$\beta_{15}^{3}\pm\beta_{j}^{4}, j =2,3,4,5,6,7,11,12;$ $\beta_{i}^{1}\pm\beta_{1}^{6}, j =1,2,\ldots,15;$
$\beta_{1}^{2}\pm\beta_{i}^{4}, j =1,2,\ldots,15;$ $\beta_{i}^{3}\pm\beta_{1}^{5}, j =1,2,\ldots,15$.
\smallskip

By applying the conclusion of Theorem \ref{BasicTheorem} we can find all subspaces on which the vectors are structural equigeodesics.
Some of them are the following:

\smallskip

\begin{center}
$
\left\{\begin{matrix}
\fr{U}_{\beta_{1}^1}\oplus_{j=2}^{5}\fr{U}_{\beta_{j}^{3}}\oplus_{k=12}^{15}\fr{U}_{\beta_{k}^{3}}  & \fr{U}_{\beta_1^1}\oplus_{j=1}^{4}\fr{U}_{\beta_{j}^4}\oplus_{k=12}^{15}\fr{U}_{\beta_{k}^4}  & 
\fr{U}_{\beta_{1}^6}\oplus_{i=1}^{15}\fr{U}_{\beta_{i}^{1}}
\\
\fr{U}_{\beta_{1}^5}\oplus_{i=1}^{15}\fr{U}_{\beta_i^3} &  \fr{U}_{\beta_{1}^3}\oplus_{i=1}^{4}\fr{U}_{\beta_i^4}\oplus_{k=12}^{15}\fr{U}_{\beta_{k}^4}  & 
\fr{U}_{\beta_{5}^3}\oplus\fr{U}_{\beta_{4}^4}\oplus\fr{U}_{\beta_{7}^4}\oplus_{i=9}^{14}\fr{U}_{\beta_i^4}   
\\
\fr{U}_{\beta_{1}^2}\oplus_{i=1}^{15}\fr{U}_{\beta_i^4} &
\fr{U}_{\beta_{5}^1}\oplus\fr{U}_{\beta_{12}^3}\oplus\fr{U}_{\beta_{12}^4} &
\fr{U}_{\beta_{1}^9}\oplus\fr{U}_{\beta_5^3}\oplus\fr{U}_{\beta_{10}^3}\oplus\fr{U}_{\beta_{15}^3}\oplus\fr{U}_{\beta_{5}^4}\oplus\fr{U}_{\beta_{10}^4}\oplus\fr{U}_{\beta_{15}^4}  
\\
 &  
\end{matrix} \right\}
$
$\cdots \cdots \cdot$

$
\left\{\begin{matrix}
\fr{U}_{\beta_{11}^1}\oplus\fr{U}_{\beta_{14}^3}\oplus\fr{U}_{\beta_{3}^4}\oplus\fr{U}_{\beta_{8}^4}\oplus\fr{U}_{\beta_{9}^4}  & \fr{U}_{\beta_{14}^3}\oplus_{i=3}^{7}\fr{U}_{\beta_i^4}\oplus\fr{U}_{\beta_{1}^4}\oplus\fr{U}_{\beta_{11}^4}\oplus\fr{U}_{\beta_{13}^{4}} 
\\
\fr{U}_{\beta_{3}^1}\oplus\fr{U}_{\beta_{4}^3}\oplus\fr{U}_{\beta_{8}^4}\oplus\fr{U}_{\beta_{11}^4}\oplus\fr{U}_{\beta_{12}^4}\oplus\fr{U}_{\beta_{15}^4} &  \fr{U}_{\beta_{14}^1}\oplus_{i=3}^{5}\fr{U}_{\beta_i^4}\oplus\fr{U}_{\beta_{1}^4}\oplus\fr{U}_{\beta_{11}^4}\oplus\fr{U}_{\beta_{8}^{4}}\oplus\fr{U}_{\beta_{9}^{4}}\oplus\fr{U}_{\beta_{13}^{4}} \\
\fr{U}_{\beta_{7}^1}\oplus\fr{U}_{\beta_{8}^3}\oplus\fr{U}_{\beta_{2}^4}\oplus\fr{U}_{\beta_{5}^4}\oplus\fr{U}_{\beta_{10}^4}\oplus\fr{U}_{\beta_{11}^4} &  \fr{U}_{\beta_{15}^1}\oplus\fr{U}_{\beta_{14}^3}\oplus\fr{U}_{\beta_{3}^4}\oplus\fr{U}_{\beta_{4}^4} \oplus\fr{U}_{\beta_{5}^4}\oplus\fr{U}_{\beta_{11}^4} 
\\
\fr{U}_{\beta_{3}^1}\oplus\fr{U}_{\beta_{3}^3}\oplus\fr{U}_{\beta_{2}^4}\oplus\fr{U}_{\beta_{5}^4}\oplus\fr{U}_{\beta_{8}^4}\oplus\fr{U}_{\beta_{11}^4}\oplus\fr{U}_{\beta_{14}^4} & 
\fr{U}_{\beta_{7}^3}\oplus\fr{U}_{\beta_{3}^1}\oplus\fr{U}_{\beta_{11}^1}\oplus\fr{U}_{\beta_{12}^1}\oplus\fr{U}_{\beta_{2}^4}
\\
\fr{U}_{\beta_{7}^3}\oplus\fr{U}_{\beta_{5}^1}\oplus\fr{U}_{\beta_{6}^1}\oplus\fr{U}_{\beta_{4}^4}\oplus\fr{U}_{\beta_{7}^4} &
\fr{U}_{\beta_{7}^3}\oplus\fr{U}_{\beta_{3}^1}\oplus\fr{U}_{\beta_{8}^1}\oplus\fr{U}_{\beta_{14}^1}\oplus\fr{U}_{\beta_{5}^4}
\\
\fr{U}_{\beta_{7}^3}\oplus\fr{U}_{\beta_{3}^1}\oplus\fr{U}_{\beta_{6}^1}\oplus\fr{U}_{\beta_{8}^4}\oplus\fr{U}_{\beta_{15}^4}  &
\fr{U}_{\beta_{7}^3}\oplus\fr{U}_{\beta_{8}^1}\oplus\fr{U}_{\beta_{12}^1}\oplus\fr{U}_{\beta_{10}^4}\oplus\fr{U}_{\beta_{15}^4}
\\
\fr{U}_{\beta_{7}^3}\oplus\fr{U}_{\beta_{6}^1}\oplus\fr{U}_{\beta_{11}^1}\oplus\fr{U}_{\beta_{7}^4}\oplus\fr{U}_{\beta_{8}^4} &
\fr{U}_{\beta_{7}^3}\oplus\fr{U}_{\beta_{11}^1}\oplus\fr{U}_{\beta_{12}^1}\oplus\fr{U}_{\beta_{14}^1}\oplus\fr{U}_{\beta_{4}^4}\oplus\fr{U}_{\beta_{8}^4}\oplus\fr{U}_{\beta_{8}^4}
\\
\fr{U}_{\beta_{1}^1}\oplus\fr{U}_{\beta_{3}^3}\oplus\fr{U}_{\beta_4^3}\oplus_{i=12}^{15}\fr{U}_{\beta_{i}^3}\oplus\fr{U}_{\beta_2^4}\oplus\fr{U}_{\beta_3^4}\oplus_{i=12}^{15}\fr{U}_{\beta_i^4} &
\fr{U}_{\beta_{12}^1}\oplus\fr{U}_{\beta_1^3}\oplus\fr{U}_{\beta_{9}^3}\oplus\fr{U}_{\beta_{15}^3}\oplus\fr{U}_{\beta_{1}^4}\oplus\fr{U}_{\beta_{9}^4}\oplus\fr{U}_{\beta_{15}^4} 
\\
\cdots & \cdots
\end{matrix} \right\}
$
\end{center}

\subsection{Structural Equigeodesic vectors on the flag manifold $\E_8/(\E_6\times\U(1)\times\U(1))$} 
Let $\Pi = \{\al_1, \al_2, \al_3, \al_4,$ $\al_5,$ $\al_6, \al_7, \al_8\}$ be a system of simple roots for $\E_8$ with highest $\widetilde{\alpha}= 2\alpha_1+3\alpha_2+4\alpha_3+5\alpha_4+6\alpha_5+4\alpha_6+2\alpha_7+3\alpha_8$. The flag manifold $\E_8/(\E_6\times\U(1)\times\U(1))$ is determined by $\Pi_K = \Pi\backslash \{\al_{1}, \al_{2}\}$.  From Table 1 we have that the positive $\fr{t}$-roots are given by  $R_{\fr{t}}^{+} = \{  \overline{\al}_1, \ \overline{\al}_1, \  \overline{\al}_1+ \overline{\al}_2,   \ \overline{\al}_1+  2\overline{\al}_2,  \  \overline{\al}_1+  3\overline{\al}_2,  \  2\overline{\al}_1+ 3 \overline{\al}_2 \}$.  According to Proposition \ref{onetone} (1), we obtain the decomposition (\ref{dec1}) where the sumbodules $\fr{m}_i$ are defined by (\ref{compon1}).  The sets $R^{\fr{n}}(j_1, j_2) = \{\sum_{i=1}^{8}c_i\al_i\in R_{M}^{+} : c_{1} = j_1, \ c_2 = j_2\}$ are given explicitly as follows:
\begin{eqnarray*}
R^{\fr{n}}(1,0) &=& \{e_1-e-2\} = \{\beta_{1}^1\}\\
R^{\fr{n}}(0,1) &=& \{e_2-e_i, i=3,4,5,6,7,8,\,  e_2+e_3+e_i, i=4,5,6,7,8, \, e_2+e_4+e_i, i=5,6,7,8, \\
&& e_2+e_5+e_i, i=6,7,8, \,  e_2+e_6+e_i, i=7,8, \, e_2+e_7+e_8, -(e_1+e_i+e_9), i=3,4,5,6,7,8\} \\
&& = \{\beta_{1}^2, \ldots, \beta_{27}^2\}\\
R^{\fr{n}}(1,1) &=& \{e_1-e_i, i=3,4,5,6,7,8, \, e_1+e_3+e_i, i=4,5,6,7,8, \,  e_1+e_4+e_i, i=5,6,7,8, \\ 
&& e_1+e_5+e_i, i=6,7,8, \, e_1+e_6+e_i, i=7,8, e_1+e_7+e_8, -(e_2+e_i+e_9), i=3,4,5,6,7,8\} \\
&& = \{\beta_{1}^{3},\ldots,\beta_{27}^3\}\\
R^{\fr{n}}(1,2) &=& \{e_i-e_9, \, e_1+e_2+e_i, i=3,4,5,6,7,8, \, -(e_3+e_i+e_9), i=4,5,6,7,8, \\
&& -(e_4+e_i+e_9), i=5,6,7,8, \, -(e_5+e_i+e_9), i=6,7,8, \, -(e_6+e_i+e_9), i=7,8, \\
&& -(e_7+e_8+e_9)\} = \{\beta_1^{4},\ldots,\beta_{27}^{4}\}\\
R^{\fr{n}}(1,3) &=& \{e_2-e_9\} = \{\beta_1^{5}\}\\
R^{\fr{n}}(2,3) &=& \{e_1-e_9\} = \{\beta_1^{6}\}
\end{eqnarray*}

Below we list some roots which satisfy Theorem \ref{BasicTheorem1}:\\
$\beta_{1}^{1}\pm\beta_{j}^{4}, \beta_{1}^{6}\pm\beta_{j}^{2}, \beta_{1}^{5}\pm\beta_{j}^{3}, j = 1,2\ldots, 27;$
$\beta_{1}^{2}\pm\beta_{j}^{3}, j =2,3,4,5,6,12,\ldots,22;$ 
$\beta_{2}^{2}\pm\beta_{j}^{3}, j =1,3,4,5,6,8,9,10,$ $11,16,17,18,19,20,21,23; $
$\beta_{3}^{2}\pm\beta_{j}^{3}, j =1,2,4,5,\ldots,11,13,14,15,19,20,21,24;$
$\beta_{4}^{2}\pm\beta_{j}^{3}, j =1,2,3,5,6,7,8,10,$ $11,12,14,15,17,18,21,25;$
$\beta_{5}^{2}\pm\beta_{j}^{3}, j =1,\ldots,9,11,12,13,$ $15,16,18,20,26;$
$\beta_{6}^{2}\pm\beta_{j}^{3}, j =1,2,3,4,5,7,8,9,10,$ $12,13,14,16,17,19,27;$
$\beta_{7}^{2}\pm\beta_{j}^{3}, j =3,4,5,6,8,9,10,11,$ $12,13,14,15,22\ldots,27,$ $\cdots$ 

\medskip

\noindent
$\beta_{1}^{2}\pm\beta_{j}^{4}, j =2,3,4,5,6,7, $ $18,\ldots,27;$
$\beta_{2}^{2}\pm\beta_{j}^{4}, j =1,3,4,5,6, 8,$ $14,15,16,17,22\ldots,27;$
$\beta_{3}^{2}\pm\beta_{j}^{4},$ $j =1,2,4,5,6,13,\ldots,17,19,20,21,25,26,$ $27;$
$\beta_{4}^{2}\pm\beta_{j}^{4},$ $j =1,2,3,5,6,10,13,14,$ $16,17,18,20,21,23,24,27;$
$\beta_{5}^{2}\pm\beta_{j}^{4},$ $j =1,2,3,4,6,11,13,14,15,17,$ $18,19,21,22,24,26;$
$\beta_{6}^{2}\pm\beta_{j}^{4},$ $j =1,2,3,4,5,12,13,14,15,16,18,18,20,$ $22,23,25;$
$\beta_{7}^{2}\pm\beta_{j}^{4},$ $j =3,4,5,6,9,$ $11,12,\ldots,18,$ $\cdots$ 

\medskip

\noindent
$\beta_{1}^{3}\pm\beta_{j}^{4}, j =2,3,4,5,6,7,18,19,20,\ldots,27;$
$\beta_{2}^{3}\pm\beta_{j}^{4}, j =1,3,4,5,6,8,14,15,16,17,22,\ldots,27;$
$\beta_{3}^{3}\pm\beta_{j}^{4}, j =1,2,4,5,6,9,13,15,16,17,19,20,21,25,26,27;$
$\beta_{4}^{3}\pm\beta_{j}^{4}, j =1,2,3,5,6,10,13,14,16,17,18,20,21,$ $23,24,27;$
$\beta_{5}^{3}\pm\beta_{j}^{4}, j =1,2,3,4,6,11,13,14,15,17,18,19,21,22,24,27;$
$\beta_{6}^{3}\pm\beta_{j}^{4}, j =1,2,3,4,5,13,14,15,16,$ $18,19,20,22,$ $23,25,27;$
$\beta_{7}^{3}\pm\beta_{j}^{4}, j =3,4,5,6,9,10,11,12,14,\ldots,21,$ $\cdots$

From these roots we can find some subspaces for which the vectors are structural equigeodesics. 
In particular we have:\\
$
\left\{\begin{matrix}
\fr{U}_{\beta_{1}^1}\oplus_{i=1}^{27}\fr{U}_{\beta_{i}^{4}} \, ,\, \fr{U}_{\beta_{1}^6}\oplus_{i=1}^{27}\fr{U}_{\beta_{i}^{2}} & \fr{U}_{\beta_{1}^5}\oplus_{i=1}^{27}\fr{U}_{\beta_{i}^{3}} 
\\
\fr{U}_{\beta_{1}^2}\oplus\fr{U}_{\beta_{2}^{3}}\oplus_{i=3}^{6}\fr{U}_{\beta_{i}^{4}}\oplus_{j=22}^{27}\fr{U}_{\beta_{j}^{4}} & \fr{U}_{\beta_{2}^2}\oplus\fr{U}_{\beta_{1}^{3}}\oplus_{i=3}^{6}\fr{U}_{\beta_{i}^{4}}\oplus_{j=22}^{27}\fr{U}_{\beta_{j}^{4}} 
\\
\fr{U}_{\beta_{2}^2}\oplus\fr{U}_{\beta_{3}^{3}}\oplus\fr{U}_{\beta_{1}^{4}}\oplus_{i=4}^{6}\fr{U}_{\beta_{i}^{4}}\oplus_{j=15}^{17}\fr{U}_{\beta_{j}^{4}}\oplus_{k=25}^{27}\fr{U}_{\beta_{k}^{4}} & \fr{U}_{\beta_{2}^6}\oplus\fr{U}_{\beta_{7}^{3}}\oplus_{i=3}^{5}\fr{U}_{\beta_{i}^{4}}\oplus_{j=14}^{16}\fr{U}_{\beta_{j}^{4}}\oplus_{k=18}^{20}\fr{U}_{\beta_{k}^{4}} 
\\
\fr{U}_{\beta_{7}^2}\oplus\fr{U}_{\beta_{6}^{3}}\oplus_{i=3}^{5}\fr{U}_{\beta_{i}^{4}}\oplus_{j=13}^{16}\fr{U}_{\beta_{j}^{4}}\oplus\fr{U}_{\beta_{18}^4} &
\fr{U}_{\beta_{1}^2}\oplus_{i=2}^{6}\fr{U}_{\beta_{i}^{3}}\oplus_{j=12}^{22}\fr{U}_{\beta_{j}^{3}}
\\
\fr{U}_{\beta_{1}^2}\oplus\fr{U}_{\beta_{2}^{2}}\oplus_{i=3}^{6}\fr{U}_{\beta_{i}^{3}}\oplus_{j=16}^{21}\fr{U}_{\beta_{j}^{3}}  & 
\fr{U}_{\beta_{1}^2}\oplus\fr{U}_{\beta_{2}^{2}}\oplus\fr{U}_{\beta_{3}^{2}}\oplus_{i=4}^{6}\fr{U}_{\beta_{i}^{3}}\oplus_{j=19}^{21}\fr{U}_{\beta_{j}^{3}}
\\
\fr{U}_{\beta_{1}^2}\oplus_{i=2}^{7}\fr{U}_{\beta_{i}^{4}}\oplus_{j=18}^{27}\fr{U}_{\beta_{j}^{4}} &
\fr{U}_{\beta_{1}^2}\oplus\fr{U}_{\beta_{2}^{2}}\oplus_{i=3}^{6}\fr{U}_{\beta_{i}^{4}}\oplus_{j=22}^{27}\fr{U}_{\beta_{j}^{4}}
\\
\fr{U}_{\beta_{1}^2}\oplus\fr{U}_{\beta_{2}^{2}}\oplus\fr{U}_{\beta_{3}^{2}}\oplus_{i=4}^{6}\fr{U}_{\beta_{i}^{4}}\oplus_{j=25}^{27}\fr{U}_{\beta_{j}^{4}}  &
\fr{U}_{\beta_{1}^3}\oplus_{i=2}^{7}\fr{U}_{\beta_{i}^{4}}\oplus_{j=18}^{27}\fr{U}_{\beta_{j}^{4}}
\\
\fr{U}_{\beta_{1}^3}\oplus\fr{U}_{\beta_{2}^{3}}\oplus_{i=3}^{6}\fr{U}_{\beta_{i}^{4}}\oplus_{j=22}^{27}\fr{U}_{\beta_{j}^{4}}  & 
\fr{U}_{\beta_{1}^3}\oplus\fr{U}_{\beta_{2}^{3}}\oplus\fr{U}_{\beta_{3}^{3}}\oplus_{i=4}^{6}\fr{U}_{\beta_{i}^{4}}\oplus_{j=25}^{27}\fr{U}_{\beta_{j}^{4}}
\\
\fr{U}_{\beta_{1}^2}\oplus\fr{U}_{\beta_{2}^{2}}\oplus\fr{\beta_{4}^{3}}\oplus\fr{U}_{\beta_{3}^{4}}\oplus_{i=5}^{6}\fr{U}_{\beta_{i}^{4}}\oplus_{j=23}^{24}\fr{U}_{\beta_{j}^{4}}\oplus\fr{U}_{\beta_{27}^4} &
\fr{U}_{\beta_{1}^3}\oplus\fr{U}_{\beta_{2}^{3}}\oplus\fr{\beta_{3}^{2}}\oplus_{i=4}^{6}\fr{U}_{\beta_{i}^{4}}\oplus_{j=25}^{27}\fr{U}_{\beta_{j}^{4}}
\\
\fr{U}_{\beta_{3}^2}\oplus\fr{U}_{\beta_{4}^{2}}\oplus\fr{U}_{\beta_{5}^{3}}\oplus_{i=1}^{2}\fr{U}_{\beta_{i}^{4}}\oplus_{j=13}^{14}\fr{U}_{\beta_{j}^{4}}\oplus\fr{U}_{\beta_{6}^4}\oplus\fr{U}_{\beta_{17}^4}\oplus\fr{U}_{\beta_{21}^4} &
\fr{U}_{\beta_{4}^3}\oplus\fr{U}_{\beta_{5}^{3}}\oplus\fr{\beta_{6}^{2}}\oplus_{i=1}^{3}\fr{U}_{\beta_{i}^{4}}\oplus_{j=13}^{14}\fr{U}_{\beta_{j}^{4}}\oplus\fr{U}_{\beta_{18}^4}
\\
\cdots & \cdots
\end{matrix} \right\}
$

\end{document}